\documentclass[hidelinks,onefignum,onetabnum]{siamart250211}

\usepackage{lipsum}
\usepackage{amsfonts,amssymb}
\usepackage{graphicx}
\usepackage{epstopdf}
\ifpdf
  \DeclareGraphicsExtensions{.eps,.pdf,.png,.jpg}
\else
  \DeclareGraphicsExtensions{.eps}
\fi
\usepackage[margin=1in]{geometry}

\newsiamremark{remark}{Remark}
\newsiamremark{hypothesis}{Hypothesis}
\crefname{hypothesis}{Hypothesis}{Hypotheses}
\newsiamthm{claim}{Claim}
\newsiamremark{fact}{Fact}
\crefname{fact}{Fact}{Facts}

\headers{{lrAA}: Low-Rank Anderson Acceleration}{D. Appel\"{o} and Y. Cheng}

\title{{lrAA}: Low-Rank Anderson Acceleration\thanks{Submitted to the editors \today \funding{Research supported by DOE Office of Advanced Scientific Computing Research under the Advanced Research in Quantum Computing program, Award Number DE-SC0025424, NSF DMS-2208164, DOE grant DE-SC0023164, and Virginia Tech.}}}

\author{Daniel Appel\"{o} \thanks{Department of Mathematics, Virginia Tech, Blacksburg, VA 24061 U.S.A. (\email{appelo@vt.edu}).} 
\and Yingda Cheng \thanks{Department of Mathematics, Virginia Tech, Blacksburg, VA 24061 U.S.A. (\email{yingda@vt.edu}).}}

\usepackage{comment}

\usepackage[compatible]{algpseudocode} 
\renewcommand{\algorithmiccomment}[1]{\bgroup\hfill\small\#~#1\egroup}

\graphicspath{{figs/}}

\begin{document}
\maketitle

\begin{abstract} 
This paper proposes a new framework for computing low-rank solutions to nonlinear matrix equations arising from spatial discretization of nonlinear partial differential equations: low-rank Anderson acceleration (lrAA). lrAA is an adaptation of Anderson acceleration (AA), a well-known approach for solving  nonlinear fixed point problems, to the low-rank format. 
In particular, lrAA  carries out all linear and nonlinear operations in low-rank form with rank truncation using an adaptive truncation tolerance. We propose a simple scheduling strategy  to  update the truncation tolerance throughout the iteration according to a residual indicator. This controls the intermediate rank and iteration number effectively. To perform rank truncation for nonlinear functions, we propose a new cross approximation, which we call  Cross-DEIM, with adaptive error control that is based on the discrete empirical interpolation method (DEIM). Cross-DEIM employs an iterative update between the approximate singular value decomposition (SVD) and cross approximation. It naturally incorporates a  warm-start strategy for each lrAA iterate. We demonstrate the superior performance of lrAA applied to a range of linear and nonlinear problems, including those arising from finite difference discretizations of Laplace's equation, the Bratu problem, the elliptic Monge-Amp\'ere equation  and the Allen-Cahn equation. 
\end{abstract}



\section{Introduction}
In this paper, we propose   low-rank Anderson acceleration (lrAA) for solving nonlinear matrix equations  
\begin{equation}
\label{eq:nle}
G(X)=X, \quad X \in \mathbb{R}^{m \times n}.
\end{equation}
In \eqref{eq:nle}, the function $G(\cdot)$ is a nonlinear function that, for example, can come from spatial discretization of a nonlinear PDE. We assume that an element of $G(X):$ $G(X)(i,j)$ is a nonlinear function of $X(i,j)$ and its nearby neighbors, i.e. the underlying discretization uses a local stencil.

In certain applications (e.g. problems with diffusion), the analytic solution is approximately low-rank, and can be approximated by a low-rank matrix with rank $r \ll O(m, n).$ 
A rank $r$ matrix of size $m \times n$ can be represented by $(m+n-r)r$ degrees of freedom which is much less than the $mn$ individual matrix elements. Thus, the goal of low-rank methods is to directly obtain a low-rank approximation with sublinear storage and computational cost, i.e. we design methods that have cost and storage proportional to $O(m+n)$ instead of $O(mn).$ In dimensions three and higher, the low-rank matrix concept is replaced by a low-rank tensor   in some compressed tensor format like tensor train (TT) or Tucker tensor.  Thus, the idea in this work can be generalized to higher dimensions. 

Developing low-rank computational methods for matrix equations is challenging. Even when $G$ is linear \cite{simoncini2016computational}, the ``best" method is chosen case-by-case depending on the specific form of the equation. For nonlinear equations, low-rank approximation of PDE is largely under-explored. \cite{sutti2024implicit} proposed a Riemannian optimization based fixed rank method. In \cite{rodgers2023implicit,adak2024tensor}, the authors proposed to use Newton method with TT-GMRES to solve the nonlinear equation when the unknown is represented in a TT format. This approach uses Newton method with  low-rank Krylov method \cite{kressner2011low,palitta2021convergence, simoncini2023analysis}, that  uses low-rank truncation within a standard Krylov method such as GMRES.  Another approach, \cite{naderi2024cur}, uses a sparse residual collocation scheme that requires the residual to vanish at certain selected columns and rows of the matrix.  

Here we take the approach by considering  the fixed-point iteration
\begin{equation} \label{eq:Picard}
X_{k+1} = G(X_k).
\end{equation}
To accelerate   the fixed point iteration, we use  Anderson acceleration (AA) \cite{anderson1965iterative, saad2024acceleration}.   AA can be used to speed up the convergence of Picard iteration, or even calculate a fixed point when the Picard iterates diverge \cite{pollock2019anderson}.  While Picard iteration only uses the current iterate to calculate the next one, $X_{k+1}$ in AA is a weighted sum of the previous $\min(k,\hat{m})+1$ iterates and residuals, where $\hat{m}$ is the window size of AA. The weighted sum is chosen so that it minimizes a linearized residual~\cite{anderson1965iterative} in the next iteration. 
It is closely related to Pulay mixing~\cite{pulay1980convergence} and DIIS (direct inversion on the iterative subspace)~\cite{kudin2002black,rohwedder2011analysis}. AA has been studied  in terms of convergence and efficient formulations \cite{walker2011anderson,toth2015convergence,evans2018proof,zhang2018globally,pollock2019anderson,de2022linear,rebholz2023effect,de2024anderson}.  AA has close connection to GMRES   when the fixed-point operator is linear~\cite{walker2011anderson}.  AA can be viewed as a multisecant quasi-Newton method~\cite{fang2009two,lin2013elliptic} and is also related to traditional series acceleration methods~\cite{weniger2001nonlinear}. 

We adapt AA to the low-rank format, which means we carry out all the linear and nonlinear operations   in low-rank form. Specifically, the unknown and the intermediate iterates are all stored in their SVD form. All linear and nonlinear operations are followed by  rank truncation to control the rank. This truncation (sometimes called rounding) is the main idea behind low-rank solvers for linear problems: low-rank Krylov methods \cite{tobler2012low, simoncini2023analysis} and iterative thresholding algorithms \cite{bachmayr2017iterative}.
The thrust of research on low-rank Krylov methods or any low-rank iterative methods is how to control the rapid growth of intermediate rank by truncation.   Frequent restart and a good  preconditioner are also deemed essential \cite{bachmayr2023low,meng2024preconditioning}. 

The AA framework provides a natural rank limiting strategy by using a finite window size. Moreover, to further control intermediate rank inflation in lrAA, we propose a simple scheduling strategy that adaptively choose the truncation tolerance throughout the iteration according to a residual indicator.
We demonstrate, computationally, that scheduling can significantly reduce computational cost by limiting the rank and lowering the number of iterations to converge and results in an almost monotonically increasing behavior of the intermediate solution rank.   

One of the main challenges for efficient low-rank methods for (\ref{eq:nle}) arises from nonlinearity. To achieve sublinear computational scaling, we use cross approximation \cite{goreinov1997theory,goreinov2001maximal} for computing low-rank approximations to terms like $G(X)$. We develop a new cross approximation, called Cross-DEIM, which performs   iterative updates between the approximate singular value decomposition (SVD) and cross approximation. Cross-DEIM is an adaptive cross approximation which gives an approximation to  matrix SVD based on a given tolerance. We verify the performance of Cross-DEIM by  testing it on benchmark problems for matrix approximation and parametric matrix approximation. For a fixed point iteration approaching convergence we expect the consecutive iterates to be ``similar". Cross-DEIM thus allows us to use the previous iterate for warm-start, which is computationally very efficient, especially when coupled with lrAA.   

The rest of the paper is organized as follows: in Section \ref{sec:lraa}, we describe the lrAA method. In Section \ref{sec:cd}, we give the details of the Cross-DEIM algorithm. Section \ref{sec:num} contains numerical tests for Cross-DEIM and lrAA methods. Section \ref{sec:conclude} concludes the paper.

\section{Low-rank Anderson acceleration}
\label{sec:lraa}

\subsection{Review of Anderson acceleration}
The convergence of the Picard iteration (\ref{eq:Picard}) is only guaranteed when certain assumptions hold on $G$ as well as the initial iterate $X_0$, and even then, its convergence is typically linear~\cite[Chap.~4.2]{kelley1995iterative}. To promote faster convergence, AA computes $X_{k+1}$ using the previous $\min(k,\hat{m}) + 1$ iterates and residuals, where $\hat{m}$ is the chosen window size. 

In the original AA algorithm, a constrained optimization problem is solved in each iteration. It is also possible to formulate AA as an unconstrained optimization problem, see \cite{walker2011anderson}. It is this formulation we use here. The   algorithm for a vector valued fixed point iteration $x_{k+1} = g(x_k) \in \mathbb{R}^n$ is presented in Algorithm~\ref{alg:AA2}. (For simplicity of presentation, we assume the relaxation parameter $\beta=1$ in AA.) 

\begin{algorithm}[ht]
   \caption{Unconstrained variant of Anderson acceleration in $\mathbb{R}^n$.} \label{alg:AA2}
  \begin{algorithmic}
  \STATE {\bf Input:}   $x_0\in \mathbb{R}^n$, window size $\hat{m} \ge 1.$
   \STATE {\bf Output:} $x_k \in \mathbb{R}^n$ as an approximate solution to $x=g(x)$.
  \STATE $x_1 = g(x_0)$.
\FOR{$k = 1, 2, \ldots$ until convergence}
       \STATE $\hat{m}_k = \min(\hat{m},k)$.
    \STATE Set $D_k = (\Delta f_{k-\hat{m}_k},\ldots,\Delta f_{k-1})$, where $\Delta f_i = f_{i+1}-f_i$ and $f_i$ = $g(x_i) - x_i$.
\STATE Solve $
\gamma^{(k)} = {\rm arg}\!\!\min_{v\in\mathbb{R}^{\hat{m}_k}} \| D_k v - f_k \|, \qquad \gamma^{(k)} = (\gamma^{(k)}_0,\ldots,\gamma^{(k)}_{\hat{m}_k-1})^T.$ 
\STATE 
$x_{k+1}   = g(x_k) -  \sum_{i=0}^{\hat{m}_k-1}  \gamma_i^{(k)} \left[g(x_{k-\hat{m}_k+i+1}) - g(x_{k-\hat{m}_k+i}) \right].$
\ENDFOR
    \end{algorithmic}
      \label{algo:AA}
\end{algorithm}

\subsection{Low-rank Anderson acceleration}
In this subsection, we propose the lrAA method. The proposed method is a low-rank adaptation of the standard AA algorithm, Algorithm \ref{algo:AA}, and is defined in Algorithm \ref{algo:AAlr}. The algorithm will output an approximate solution to the fixed point problem \eqref{eq:nle} in its  SVD form:  $X_k=U S V^T$, with $U\in \mathbb{R}^{m \times r},$ $S \in \mathbb{R}^{r \times r}$ being diagonal and $V\in \mathbb{R}^{r \times n}$. $r\le r_{\rm max},$ and the output meets the   tolerance in the sense that $\|G(X_k)-X_k\| \le {\rm TOL}$ as specified by the input parameter ${\rm TOL}$. Here, we note that in this paper we always take the matrix Frobenius norm and denote it by $\|\cdot\|.$ We now describe the different elements of lrAA.

\begin{algorithm}[ht!]
  \caption{lrAA for nonlinear matrix equation $G(X)=X.$     \label{alg:lrAA} }
  \begin{algorithmic}[1]
  \STATE {\bf Input:} $X_0=U_0 S_0 (V_0)^T$, window size $\hat{m} \ge 1,$  scheduling parameter $\theta \in(0,1)$, tolerance ${\rm TOL},$ max rank  $r_{\rm max},$ truncation parameter $\epsilon_F.$ 
     \STATE {\bf Output:} Approximate solution $X_k$ to the fixed point problem $G(X)=X$ in its SVD form.
  \STATE $\epsilon_{G} = 10^{-2}$ \COMMENT{Choose $\epsilon_{G}$ so that $G_0$ has low rank.}

\STATE $X_1 = G_0=\textbf{\textrm{Cross-DEIM}}(G(X_0),U_0,V_0,\epsilon_{G},r_{\rm max})$. 
\STATE   $\rho_0 = \| G_0 - X_0 \|$. \COMMENT{approximate residual as difference of low-rank matrices}
\FOR{$k = 1, 2, \ldots$  }
\STATE 
$
G_k=\textbf{\textrm{Cross-DEIM}}(G(X_k),U_k,V_k,\epsilon_G,r_{\rm max}).$
\STATE   $\rho_k = \| G_k -X_k\|$. \COMMENT{approximate residual as difference of low-rank matrices}
\STATE  $\hat{m}_k = \min(\hat{m},k)$.
\STATE Set $D_k = (\Delta F_{k-\hat{m}_k},\ldots,\Delta F_{k-1})$, where $\Delta F_i =  \mathcal{T}^{\rm round}_{\epsilon_F,r_{\rm max}}(F_{i+1}-F_i$) and \mbox{$F_i =  \mathcal{T}^{\rm round}_{\epsilon_F,r_{\rm max}}(G_i - X_i)$}.
\STATE Solve 
\begin{equation*}
\gamma^{(k)} = {\rm arg}\!\!\min_{v\in\mathbb{R}^{\hat{m}_k}}\|f_k- D_k v\|, \qquad \gamma^{(k)} = (\gamma^{(k)}_0,\ldots,\gamma^{(k)}_{\hat{m}_k-1})^T. \label{eq:distance1}
\end{equation*}
\STATE  $X_{k+1}  = \textbf{\textrm{Cross-DEIM}}(G_k -  \sum_{i=0}^{\hat{m}_k-1}  \gamma_i^{(k)} \left[G_{k-\hat{m}_k+i+1} - G_{k-\hat{m}_k+i} \right] ,U_k,V_k,\epsilon_G,r_{\rm max}).$
\STATE Set $\epsilon_G = \theta \rho_k$. \COMMENT{scheduling: update the truncation tolerance} 
\IF {$\rho_k < {\rm TOL}$} 
\STATE Exit and return $X_{k+1}$.
\ENDIF
\ENDFOR
    \end{algorithmic}
    \label{algo:AAlr}
\end{algorithm}

First we note that throughout the algorithm quantities like $X_k$ and $G_k$ should be interpreted as their SVD representation, we never form the $m \times n$ matrices, but instead we always operate on their factors. A core component of the algorithm is the Cross-DEIM method, which will be described in Section \ref{sec:cd} and is defined in Algorithm \ref{algo:crossDEIM}. What is important here is that it returns the SVD representation of a low-rank matrix of a nonlinear function (as in lines 4 and 7) or a linear combination of low-rank matrices (as in line 12) according to a prescribed tolerance $\epsilon_G$ and max rank $r_{\rm max}$. (Here for line 12, we use Cross-DEIM as an approximate rounding operation with better efficiency than the standard rounding method: Algorithm \ref{algo:round} for large window size $\hat{m}$, but rounding can also be used as an alternative.) In Cross-DEIM, we use  the singular vectors $U_k, V_k$ from the previous iterate $X_k$ to warm-start the cross approximations of both $G(X_k)$ and $X_{k+1}.$  This works very well in practice and the mechanism will be explained in details in Section \ref{sec:cd}.

The parameter $\epsilon_G$ is a  tolerance such controlling the error in the   Cross-DEIM approximation. It serves a similar role as the  truncation tolerance in a  truncated SVD. The quantity $\rho_k$ is a residual used to determine when to stop the iteration. This quantity is the norm of the difference of two low-rank matrices, and can be efficiently calculated using similar spirit of Algorithm \ref{algo:round}. The residual $\rho_k$ is  used to adaptively schedule the tolerance according to the user defined  scheduling parameter $\theta,$ which is critical to the good performance of lrAA. We propose the following simple strategy to take $\epsilon_G$ to be proportionate to $\rho_k.$ Precisely, given $\theta \in (0,1)$ and $\rho_k,$ we set $\epsilon_G=\theta \rho_k$ (line 13). As we will show in the numerical experiments section, this simple strategy works well to keep the intermediate rank low. In the computations in this paper we take $\theta \approx 0.1-0.5$, this appears to work well for the examples considered here  (we   note that more advanced strategies from \cite{bachmayr2017iterative} could be extended to lrAA but we leave this as future work). The initial value for $\epsilon_G$ should be chosen so that the initial iterate $X_1$ has low rank. We typically set it to $10^{-2}$. In practice, $\epsilon_G$   will gradually decrease with $\rho_k$ till $\rho_k$ meets the prescribed tolerance. Thus, qualitatively, the numerical rank of the iterates $X_k$ will gradually increase over the iteration till it reaches the  rank needed. Such  almost monotonically-increasing  rank property is highly desirable for low-rank methods from the perspective of computational efficiency. We   mention that Algorithm \ref{algo:AAlr} can also be turned into a fixed-rank version by using the $r_{\rm max}$ option. For brevity, in this paper, we focus on rank-adaptive version   and do not experiment with fixed-rank lrAA. Therefore, in our numerical experiment, we always disable $r_{\rm max}$ by setting it to be the $\min(m,n).$

Lines 9 to 11 find the solution to the least squares minimization problem using Algorithm \ref{algo:leastsquareLR} and updates the solution as a linear combination of previous $G_k$. The rounding of sums of low-rank matrices used in line 9 is described in Algorithm \ref{algo:round}. In order for the least squares solution to be accurate, we have found that the rounding tolerance  $\epsilon_F$ should be chosen quite small. In all the examples in this paper, we take it to be $10^{-12}$. Taking $\epsilon_F$ too large typically results in ill-conditioning and an increase in the number of iterations to reach ${\rm TOL}$. As for the window size, the standard AA methods recommend  $\hat{m} \approx 3-10$ \cite{kelley2022solving}. When the window size is too small, it will negatively affect the iteration number. For lrAA, as a low-rank method, when the window size is too large,  it will negatively affect both the iteration number and the efficiency of the least square solve. Therefore, based on the numerical experiments, we recommend a slightly smaller window size of about $3-5.$ 

As can be seen, besides Cross-DEIM, which will be discussed in the next section, all other elementary operations in lrAA are performed in low-rank format with operation cost on par of $O((m+n)r^2)$ if 
$r \ll O(m, n).$  Moreover, as we will show in the numerical experiment section, the iteration number of lrAA is shown to be less than their full rank counterpart (the standard AA) due to scheduling. Therefore, we expect significant gain in  computational cost and storage for lrAA compared to standard AA when the low-rank assumption holds. Finally, let's mention that in this work we consider the most basic versions of AA, i.e. we did not consider QR-updating techniques \cite{golub2013matrix} and other late improvements to AA \cite{walker2011anderson, lupo2024anderson} which can improve computational efficiency and potentially alleviate the restriction on $\epsilon_F$. We believe those techniques can be useful and will be explored in the future.   

 \begin{algorithm}[htb]
 \caption{Rounding of sum of low rank matrices, i.e.   $US V^T = \mathcal{T}^{\rm round}_{\epsilon,r_{\rm max}}(\sum_{j=1}^d U_jS_jV^T_j)$}
    \begin{algorithmic}
      \STATE {\bf Input:} low rank matrices in the form $U_jS_jV^T_j, j =1,\dots,s$,   tolerance $\epsilon,$ max rank $r_{\rm max}$
      \STATE {\bf Output:} $U, S, V$
   \STATE Let $U=[U_1,\dots ,U_s]$, $S={\sf diag}(S_1,\dots,S_s)$, $V=[V_1,\dots,V_s]$
    \STATE Perform column pivoted QR: $[Q_1,R_1,\Pi_1] = {\sf qr}(U)$, $[Q_2,R_2,\Pi_2] = {\sf qr}(V)$
   \STATE Compute the truncated SVD for the small matrix with tolerance $\epsilon$ and max rank $r_{\rm max}$: $\mathcal{T}_{\epsilon,r_{\rm max}}(R_1\Pi_1S\Pi_2^TR_2^T) = US V^T$
    \STATE  $U \leftarrow Q_1 U$, $V \leftarrow Q_2 V$
 	\end{algorithmic} 
	\label{algo:round}
\end{algorithm} 

 \begin{algorithm}[htb]
 \caption{Computing the least squares solution minimizing $\| \sum_{j=1}^{s} \gamma_j U_jS_j V_j^T - U_B S_B V_B^T \|$	\label{algo:leastsquareLR}}
    \begin{algorithmic}
      \STATE {\bf Input:} low rank matrices in the form $U_jS_jV^T_j, j =1,\dots,s$, right hand side $U_B S_B V_B^T$
            \STATE {\bf Output:} $\gamma_j, j=1, \ldots s$
      \STATE Let $U=[U_1,\dots ,U_s]$, $V=[V_1,\dots,V_s]$
    \STATE Perform column pivoted QR: $[Q_1,R_1,\Pi_1] = {\sf qr}(U)$, $[Q_2,R_2,\Pi_2] = {\sf qr}(V)$
    \STATE Set $b = {\sf vec}(Q_1^T U_B S_B V_B^T Q_2)$.
    \STATE Find the least squares $\gamma$ that minimizes the small problem $\|A \gamma - b\|$ where the $k$th column of $A$ is $a_k = {\sf vec}(R_1\Pi_1^T D_k \Pi_2 R_2^T)$,  and $D_k = {\sf diag}(0,\ldots,0,S_k,0,\ldots,0)$.  
 	\end{algorithmic} 
\end{algorithm}

\section{Cross-DEIM approximation}
\label{sec:cd}
In this section, we describe how we find a low-rank SVD approximation to a matrix $G(i,j)$, $1 \leq i \leq m, 1 \leq j \leq n$ using   cross approximation by Cross-DEIM.  Low-rank approximation using matrix skeleton, pseudoskeleton, CUR factorization and cross approximation is a well-studied subject in numerical linear algebra \cite{goreinov1997theory,goreinov2001maximal,mahoney2009cur}. Given a low-rank matrix $G,$    we approximate  $G$ by selected columns (indexed by $\mathcal{J}$) and rows (indexed by $\mathcal{I}$), i.e. $$G \approx G(:,\mathcal{J}) U G(\mathcal{I},:).$$
If $U= G(:\mathcal{J})^+ G G(\mathcal{I},:)^+,$ then the resulting approximation is usually called the CUR decomposition, and it is the best approximant to $G$ in the Frobenious norm, i.e. $G(:\mathcal{J})^+ G G(\mathcal{I},:)^+={\rm arg} \min_U \|G-G(:\mathcal{J}) U G(\mathcal{I},:)\|_F$. However, this approach uses the full matrix $G$ and does therefore not satisfy the sublinear scaling. In this work, we instead use cross approximation $U=G(\mathcal{I},\mathcal{J})^+$ in its stabilized version to achieve sublinear scaling.  

The quality of cross approximation critically hinges on the choice of the column and row indices, $\mathcal{I}, \mathcal{J}.$ Popular methods for choosing  $\mathcal{I}$ and $ \mathcal{J}$ are randomized algorithms \cite{halko2011finding,chiu2013sublinear,cortinovis2025sublinear}, max volume \cite{goreinov2001maximal}, leverage score \cite{mahoney2009cur},
and DEIM (discrete empirical interpolation method) and its variant QDEIM (potentially with oversampling) \cite{sorensen2016deim,drmac2016new,peherstorfer2020stability}.  The DEIM based index selection for CUR decomposition, proposed in \cite{sorensen2016deim}, works well when the leading singular vectors of $G$ are available and the performance of DEIM type methods have been tested and shown to be advantageous compared with leverage score based index selection. Due to the availability of approximate singular vector information, DEIM index selection has recently been adopted in for low-rank methods for nonlinear PDEs and stochastic PDEs \cite{donello2023oblique,ghahremani2024deim,dektor2024collocation,ghahremani2024cross,naderi2024cur}. 
 
In this work, we focus on cross approximation with a warm-start strategy for index selection. Warm-start strategy has been recently explored in 
\cite{park2024low}, where the index set is directly recycled for parameter-dependent problems, and 
in \cite{donello2023oblique,ghahremani2024deim,dektor2024collocation,ghahremani2024cross,naderi2024cur}, where the index set from DEIM method of singular vectors from previous time step for time-dependent problem is used as the starting index set. Warm-start strategy naturally applies to lrAA and other low-rank iterative schemes (where the iteration number is the parameter). 

Our method explores the warm-start strategy by the DEIM index selection. We chose this approach because  of the readily available singular vector information from previous lrAA iterate. What distinguishes this work from \cite{donello2023oblique,ghahremani2024deim,dektor2024collocation,ghahremani2024cross,naderi2024cur} is that rather than just using the index from approximate singular vectors, we propose an iterative procedure that updates the index selection and the SVD approximation iteratively to achieve adaptive  control on the error bounds. In contrast to prior iterative index update \cite{oseledets2010tt, ghahremani2024cross} by replacement for a fixed rank output, our method contains   merging and pruning process for an adaptive process with error tolerance control.  

 The algorithm we propose  is denoted   by $\textbf{\textrm{Cross-DEIM}}(G, U_0, V_0, \epsilon,r_{\rm max},\textsf{opts}),$ where $G$ is the matrix to be approximated; $U_0, V_0$ denote singular vectors to initialize the DEIM index selection; $\epsilon$ is the error tolerance (in matrix Frobenius norm); $r_{\rm max}$ is the maximum rank; $\textsf{opts}$ are options. The output is a low-rank matrix in its SVD form.  The details of the algorithm is described in Algorithm \ref{algo:crossDEIM}.
 Algorithm \ref{algo:crossDEIM} uses the subroutines QDEIM for index selection \cite{drmac2016new} (restated in Algorithm \ref{algo:qdeim}) and a version of stabilized cross approximation \cite{donello2023oblique} stated in Algorithm \ref{algo:crossCUR}. The QDEIM subroutine takes inputs of a tall orthogonal matrix of size $k \times l$ and outputs the leading $l$ important rows by column pivoted QR. In practice, DEIM index selection \cite{sorensen2016deim} can also be used with little  qualitative difference. 

Now, we outline the main steps in the  \textrm{Cross-DEIM} algorithm. We initialize the routine with empty row and column index sets and an initial guess of approximate left and right leading singular vector matrix $U_{0}, V_{0}.$ Then at each iterate $k,$ given the current row and column index sets $\mathcal{I}_k, \mathcal{J}_k$ and approximate left and right leading singular vector matrix $U_{k-1}, V_{k-1},$ the iterative procedure updates the index sets (lines 5-19) and the SVD (line 20) until it converges. 
Specifically, the index update is done by enriching the current index sets with indices selected by QDEIM with singular vector matrices $U_{k-1}, V_{k-1}$ (lines 5-7). As a ssafeguard,we add one additional randomly chosen index if no new indices are added in the QDEIM step (lines 8-13). The stabilized cross approximation updates the singular vectors. We remove the redundant (linearly dependent) rows and columns as needed (lines 21-30). The iteration is stopped when both the difference of the consecutive updates and the low-rank indicators are less than the provided tolerance. Here, the low-rank indicator defined on lines 32-33 of Algorithm \ref{algo:crossDEIM} is based on the error estimates from \cite{donello2023oblique}. Finally, additional pruning is conducted based on the error tolerance $\epsilon$  (lines 37-38).

The cost of the algorithm will depend on $m, n,  $ the output rank $r$ as well as the iteration number and the intermediate rank throughout the iteration. In general, the main computational cost will be the call to the stabilized cross approximation in Algorithm \ref{algo:crossCUR}, which scale as $O((m+n)r^2).$ Thus, as long as the iteration number remains $O(1)$ and the max intermediate rank is a constant multiple of $r$, we expect this sublinear cost  for the overall algorithm. We demonstrate in numerical experiments in Section \ref{sec:numcd} that such assumptions are valid for matrix with rapid singular value decay. For more challenging problems with slow singular value decay, the number of iterations and the intermediate rank will be more sensitive to the choice of $U_0, V_0$ similar to other cross approximation methods. In   numerical examples from parametric matrix approximation and lrAA, we demonstrate the effectiveness of the warm-start strategy  in controlling both the iteration number and the maximum intermediate rank. Specifically, for lrAA applications tested in this paper, we observe the range of Cross-DEIM iteration number is between 2 to 4 and the max intermediate rank is a multiple smaller than two of the final rank. 

Finally, we mention that since  Cross-DEIM is a sampling method, it's possible to construct   adversarial examples for which the algorithm will fail. However, we don't expect this to happen for most applications we have in mind for lrAA, that is nonlinear PDE problems for which the target matrix $X$ demonstrates reasonably fast singular value decay. 

\begin{algorithm}[ht!]
  \caption{{\tt [U,S,V] = Cross-DEIM(G,U0,V0,$\epsilon$, $r_{\rm max}$, $\aleph_{\rm max}$, {\tt maxiter} )}\\
  Adaptive Cross-DEIM approximation to  $G \in \mathbb{R}^{m \times n}$}
  \begin{algorithmic}[1]
  \STATE {\bf Input:} Matrix $G \in \mathbb{R}^{m \times n}$, initial rank $r$ guess to the singular vector matrix $U_0 \in \mathbb{R}^{m \times r}, V_0 \in \mathbb{R}^{n \times r}$, tolerance $\epsilon$, maximum output rank $r_{\rm max}$, maximum index set cardinality $\aleph_{\rm max}$, maximum number of iterations {\tt maxiter}. 
  \STATE {\bf Output:} Approximate SVD of $G$, $U\in \mathbb{R}^{m \times r},$ $S \in \mathbb{R}^{r \times r}$, $V\in \mathbb{R}^{r \times n}$.
  \STATE Set $\mathcal{I}_0=\mathcal{J}_0 = \emptyset$.
\FOR{$k = 1, 2, \ldots$, {\tt maxiter}}
\STATE $\mathcal{I}_{k}^\ast = {\tt QDEIM}(U_{k-1})$
\STATE $\mathcal{J}_{k}^\ast = {\tt QDEIM}(V_{k-1})$ \COMMENT{ {\tt QDEIM} can be replaced by {\tt DEIM}}
\STATE $\mathcal{I}_{k} = \mathcal{I}_{k}^\ast \cup \mathcal{I}_{k-1}, \mathcal{J}_{k} \leftarrow \mathcal{J}_{k}^\ast \cup \mathcal{J}_{k-1}$ \COMMENT{Note that the index sets are ordered by {\tt QDEIM}.}
\IF[Make sure that that the index set increase by one]{$|\mathcal{I}_{k}| = |\mathcal{I}_{k-1}|$ or $k=1$} 
\STATE $\mathcal{I}_{k} = \mathcal{I}_{k}^\ast \cup \{ i_{\rm rand} \in \complement (\mathcal{I}_{k}^\ast) \}$ \COMMENT{using a random $i_{\rm rand}$ from the complement of $\mathcal{I}_{k}^\ast$.}
\ENDIF
\IF{$|\mathcal{J}_{k}| = |\mathcal{J}_{k-1}|$ or $k=1$} 
\STATE $\mathcal{J}_{k} \leftarrow \mathcal{J}_{k}^\ast \cup \{ j_{\rm rand} \in \complement (\mathcal{J}_{k}^\ast) \}$ 
\ENDIF
\IF{$|\mathcal{I}_k| > \aleph_{\rm max}$}
\STATE $\mathcal{I}_k \leftarrow \mathcal{I}_k(1:\aleph_{\rm max})$ \COMMENT{Keep the $\aleph_{\rm max}$ most important indices.} 
\ENDIF
\IF{$|\mathcal{J}_k| > \aleph_{\rm max}$}
\STATE $\mathcal{J}_k \leftarrow \mathcal{J}_k(1:\aleph_{\rm max})$ 
\ENDIF
\STATE $[U_k,S_k,V_k,r_{\rm C},r_{\rm R}] = {\tt scross}(G, \mathcal{I}_k,\mathcal{J}_k)$
\FOR{$l = 1, 2, \ldots, |\mathcal{I}_k|$}
\IF{$|(r_{\rm R})_{l}| < 10^{-12}$}
\STATE Remove element $l$ from $\mathcal{I}_k$  \COMMENT{Remove redundant rows in $R = G(\mathcal{I}_k,:)$.}
\ENDIF  
\ENDFOR 
\FOR{$l = 1, 2, \ldots, |\mathcal{I}_k|$}
\IF{$|(r_{\rm C})_{l}| < 10^{-12}$}
\STATE Remove element $l$ from $\mathcal{J}_k$ \COMMENT{Remove redundant columns in $C = G(:,\mathcal{J}_k)$.} 
\ENDIF
\ENDFOR 
\STATE $\rho = \| U_{k}S_{k}V_{k}^T- U_{k-1}S_{k-1}V_{k-1}^T \|, \ \ S_{\rm min} =  \min({\sf diag}(S_k))$
\STATE  $\eta_1 = \| (I(:,\mathcal{I}_k))^T U_k \|_2^{-1}, \ \ \eta_2 = \| V^T_{k} I(\mathcal{J}_k,:) \|_2^{-1}$ 
\IF{$\max(\rho,\min(\eta_1(1+\eta_2),\eta_2(1+\eta_1))S_{\rm min})  < \epsilon$} 
\STATE Break out of for loop  \COMMENT{Above $S_{\rm min}$ is the smallest s.v. in the $k$th approx.}
\ENDIF
\ENDFOR
\STATE Find $r^\ast$ so that $\sum_{l = r^\ast+1}^{\min(m,n)} S_l^2 < \epsilon^2$ 
\STATE Set $r = \max(\min(r^\ast,r_{\rm max}),1)$
\STATE Return $U_k(:,1:r), S_k(1:r,1:r), V_k(:,1:r)$
\end{algorithmic}
\label{algo:crossDEIM}
\end{algorithm}

\begin{algorithm}[ht!]
\caption{{\tt [$\mathcal{I}$] = QDEIM(U)} \\QDEIM index selection}
 \begin{algorithmic}[1]
 \STATE {\bf Input:}  Orthogonal matrix $U$ of size $k \times l$
\STATE {\bf Output:} Index set $\mathcal{I}$ of size $l$
\STATE  $[\sim,\sim, p]=\textrm{qr}(U^T,\textrm{'vector'})$ \COMMENT{Perform column pivoted QR on $U^T$.}
\STATE $\mathcal{I}=P(1:l)$
\STATE Return $\mathcal{I}$.
\end{algorithmic}
\label{algo:qdeim}
\end{algorithm}

\begin{algorithm}[ht!]
\caption{{\tt [U,S,V,rC,rR] = scross(G,I,J)} \\Stabilized cross  approximation of $G$}
 \begin{algorithmic}[1]
 \STATE {\bf Input:}  $G$, and two index sets $\mathcal{I}$, $\mathcal{J}$ with $k$ and $l$ elements, respectively.
\STATE {\bf Output:} Approximate SVD of $G$, $U\in \mathbb{R}^{m \times r},$ $S \in \mathbb{R}^{r \times r}$, $V\in \mathbb{R}^{r \times n}$, and two vectors to test for linear dependence $r_{\rm R}, r_{\rm C} \in \mathbb{R}^r$.
\STATE $C = G(:,\mathcal{J}) \in \mathbb{R}^{m \times k}, \ \ R = G(\mathcal{I},:) \in \mathbb{R}^{l \times n}$  
\STATE $CP_{\rm C} = QR_{\rm C}$, \ \ $R^TP_{\rm R} = ZR_{\rm R}$ \COMMENT{Perform column pivoted QR.}
\IF{$k \le l$}
\STATE Solve $Q(\mathcal{I},:) W = R$ \COMMENT{Solved using {\tt \textbackslash} if $Q(\mathcal{I},:)$ is well conditioned,}
\STATE \COMMENT{else solved by truncated SVD pseudoinverse.} 
\STATE $W = \hat{U}SV^T$  \COMMENT{Perform truncated SVD.}
\STATE $U = Q\hat{U}$
\ELSE
\STATE Solve $Z(:,\mathcal{J}) W = C^T$ \COMMENT{Solved using {\tt \textbackslash} if $Z(:,\mathcal{J})$ is well conditioned,}
\STATE  \COMMENT{else solved by truncated SVD pseudoinverse.} 
\STATE $W^T = US\hat{V}^T$  \COMMENT{Perform truncated SVD.}
\STATE $V = Z\hat{V}$
\ENDIF
\STATE Return $USV^T \approx G$ and $r_{\rm R} = {\sf diag}(P_R^T R_R P_R), r_{\rm C} = {\sf diag} (P_C^T R_C P_C)$.
\end{algorithmic}
\label{algo:crossCUR}
\end{algorithm}

\section{Numerical examples}
\label{sec:num}

We report the numerical results of Cross-DEIM and lrAA in this section. All the errors in this section refer to the errors measured in the matrix Frobenius norm. For convenience and brevity, we do not rescale the norms by the mesh-size, but note that this could be important when comparing performance of the methods under grid refinement.  

\subsection{Cross-DEIM examples}
\label{sec:numcd}

In this subsection, we test the performance of the Cross-DEIM algorithms for matrix approximation and parametric matrix approximation. We consider matrices with fast and slow singular value decay to demonstrate the behavior of the method for various scenarios. 

\subsubsection{Matrix approximation}
In this experiment we consider the approximation of $m \times n$ matrices whose elements are given by  
\begin{equation*}
G_{1}(i,j) = \frac{1}{i+j-1}, \ \ \ \  G_2(i,j) = \left(\frac{|x_i+y_j|}{2}\right)^5. 
\end{equation*}
Here $x_i = -1 + 2\frac{i-1}{m-1}$ and $y_j = -1 + 2\frac{j-1}{n-1}$. For $G_1$ we take $m=n=100$, for $G_2$ we take $m=n = 500$. $G_1$ is the Hilbert matrix which is known to have rapid singular value decays, where $G_2$ is derived from a  $C^4$  function with non-smooth feature  along $x+y=0$ on   uniform grid points. Thus, the matrix $G_2$ has slow singular value decay and poses numerical challenges for cross approximation. In this example, we benchmark the results with the truncated SVD and the adaptive cross approximation (ACA) with partial pivoting (Algorithm 4 in \cite{aca}) with the same tolerance numbers.  In the examples, we set  the Cross-DEIM parameters $r_{\rm max}, \aleph_{\rm max}, {\tt maxiter}$ to the maximally allowable values (effectively turning them off). Instead, we vary the tolerance $\epsilon$. The initial vectors $U_0, V_0$ are set as random vectors of size $m \times 1$ and $n \times 1$. To assess the overall performance of Cross-DEIM based on such random initialization, we take 100 runs, and report the mean and max ranks and errors.
 
\begin{figure}[htb]
\begin{center}
\includegraphics[width=0.35\textwidth,trim={0.0cm 0.0cm 0.0cm 0.0cm},clip]{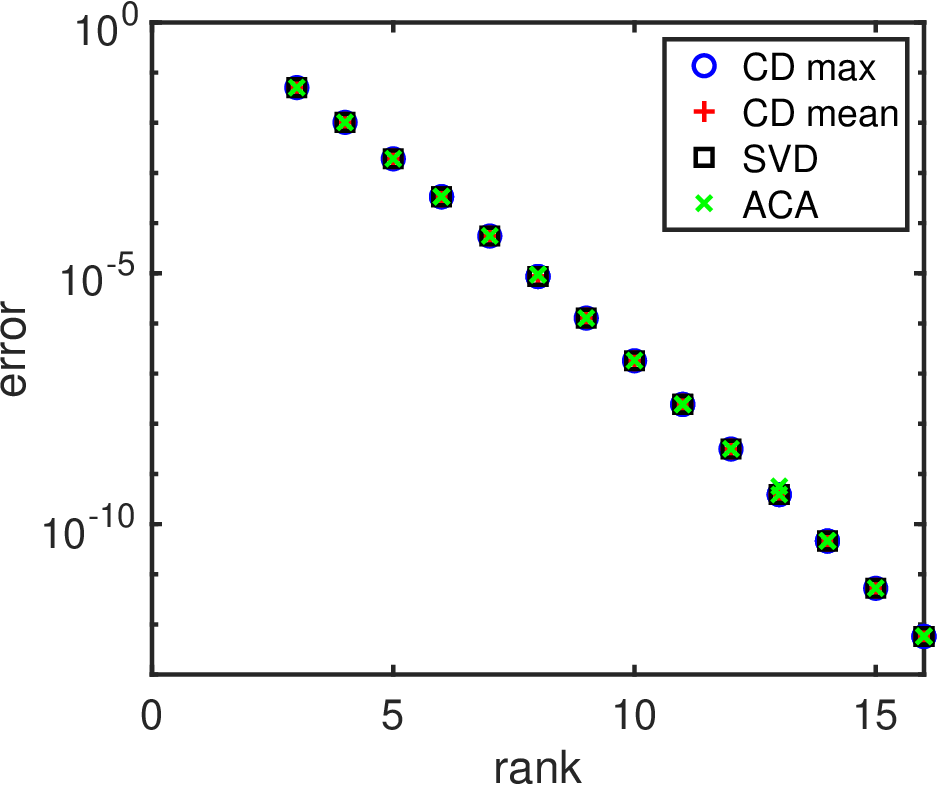}
\includegraphics[width=0.35\textwidth,trim={0.0cm 0.0cm 0.0cm 0.0cm},clip]{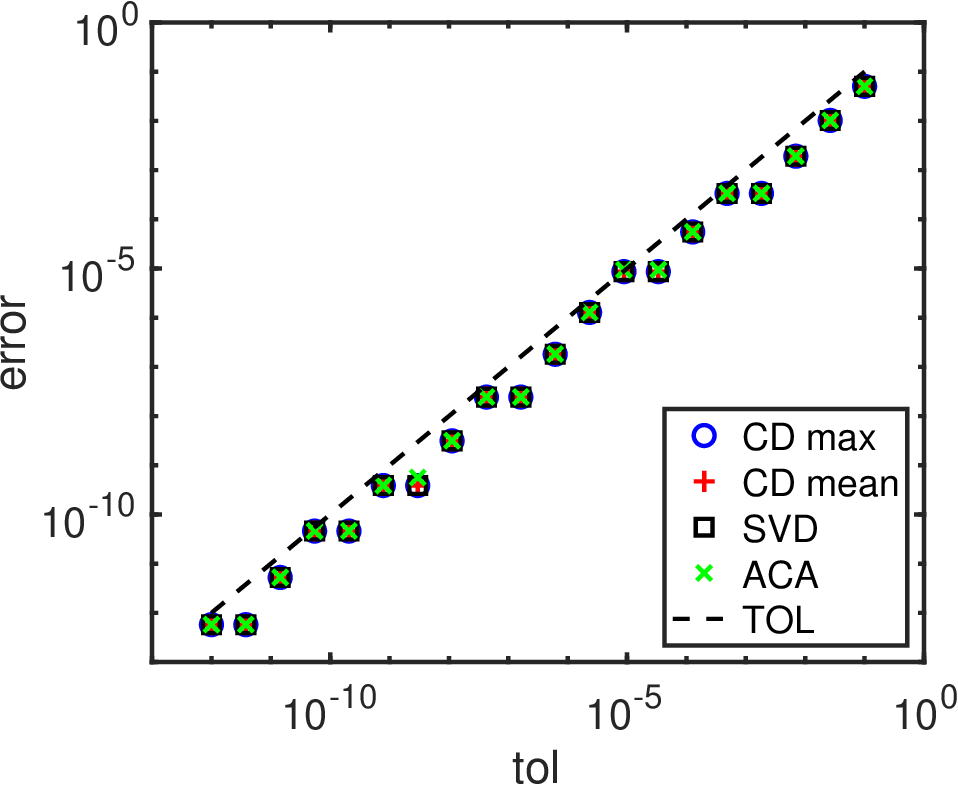}
\includegraphics[width=0.35\textwidth,trim={0.0cm 0.0cm 0.0cm 0.0cm},clip]{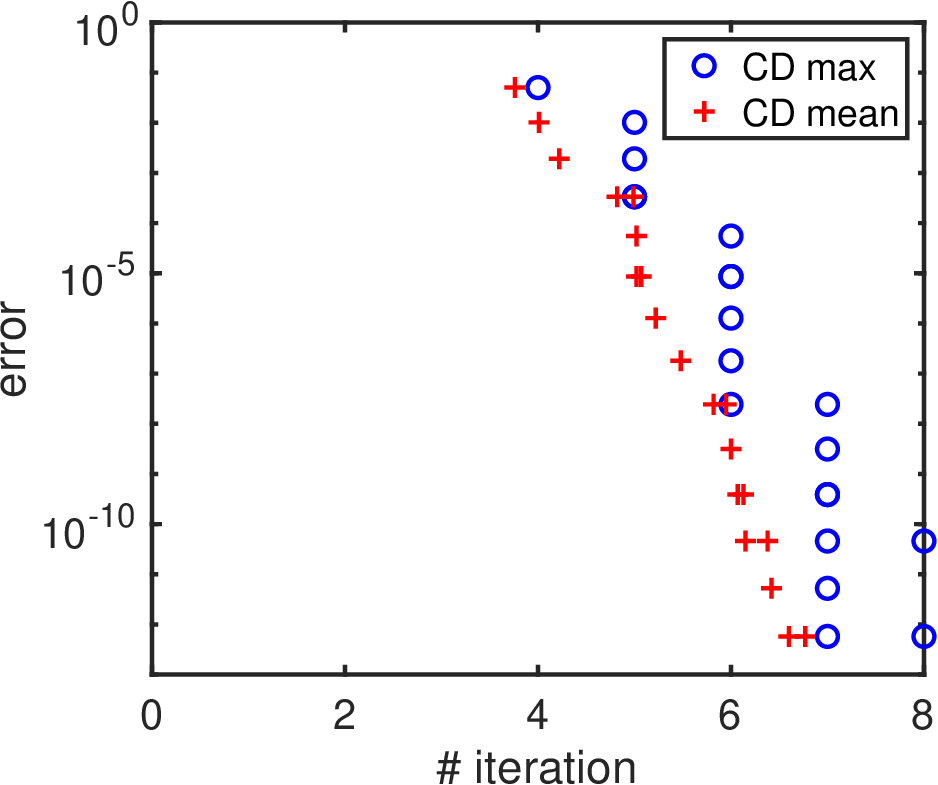}
\includegraphics[width=0.35\textwidth,trim={0.0cm 0.0cm 0.0cm 0.0cm},clip]{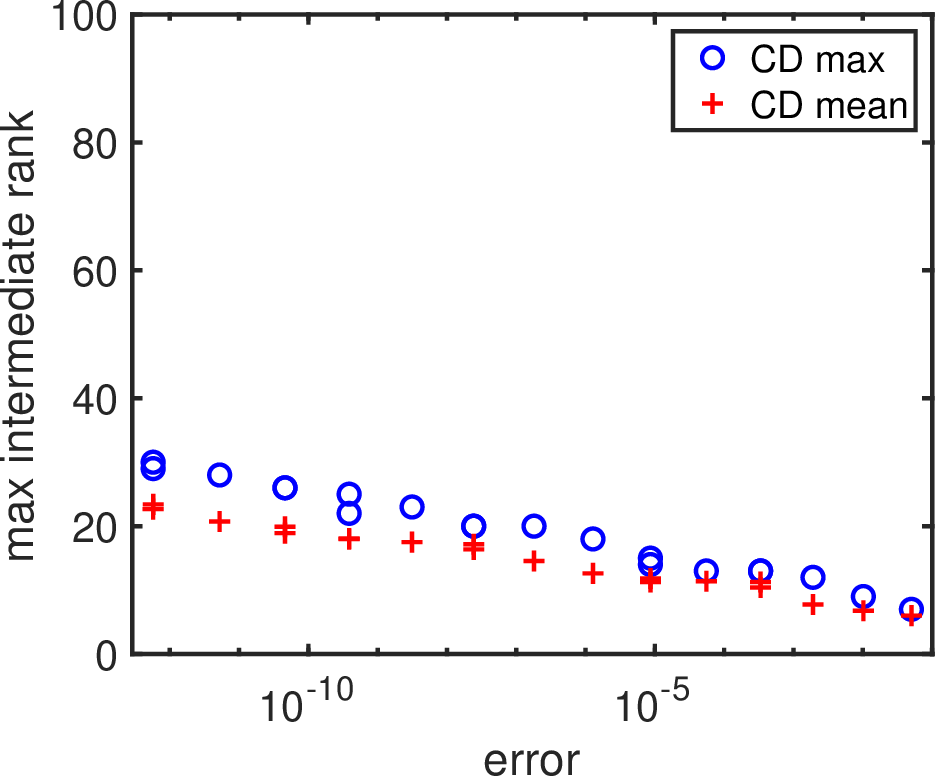}
\caption{Results for approximation of  $G_1$. In the figures, `CD' denotes Cross-DEIM. \label{fig:CD1}}
\end{center}
\end{figure}

\graphicspath{{figures/init_rank_1}}
The results for   $G_1$ with tolerance $\epsilon$ varying from $10^{-12}$ to $10^{-1}$ are reported  
in Figure \ref{fig:CD1}. In the top two subfigures, we display the  numerical   error vs   numerical rank, and the numerical error vs the requested tolerance  using Cross-DEIM, the truncated SVD and the ACA at a given tolerance. As can be seen, both Cross-DEIM and ACA demonstrate good performance, i.e. the   errors the optimal value with a given rank computed by truncated SVD (from the top left subfigure) and the numerical errors meet (are smaller than) the requested tolerance (from the top right subfigure). In particular, the results are not sensitive to the random initial vectors in the Cross-DEIM method. To benchmark the computational efficiency of the Cross-DEIM method,  in the bottom two subfigures, we plot the max and mean of iteration numbers and max intermediate rank (defined as $\max_k \max (| \mathcal{I}_k |,| \mathcal{J}_k |)$) for the Cross-DEIM method. We can see for various runs with different tolerances, the iteration numbers vary slightly from 4 to 8. Comparing the max intermediate rank with the final rank we see that for no error the ratio is greater than two. 

\begin{figure}[htb]
\begin{center}
\includegraphics[width=0.35\textwidth,trim={0.0cm 0.0cm 0.0cm 0.0cm},clip]{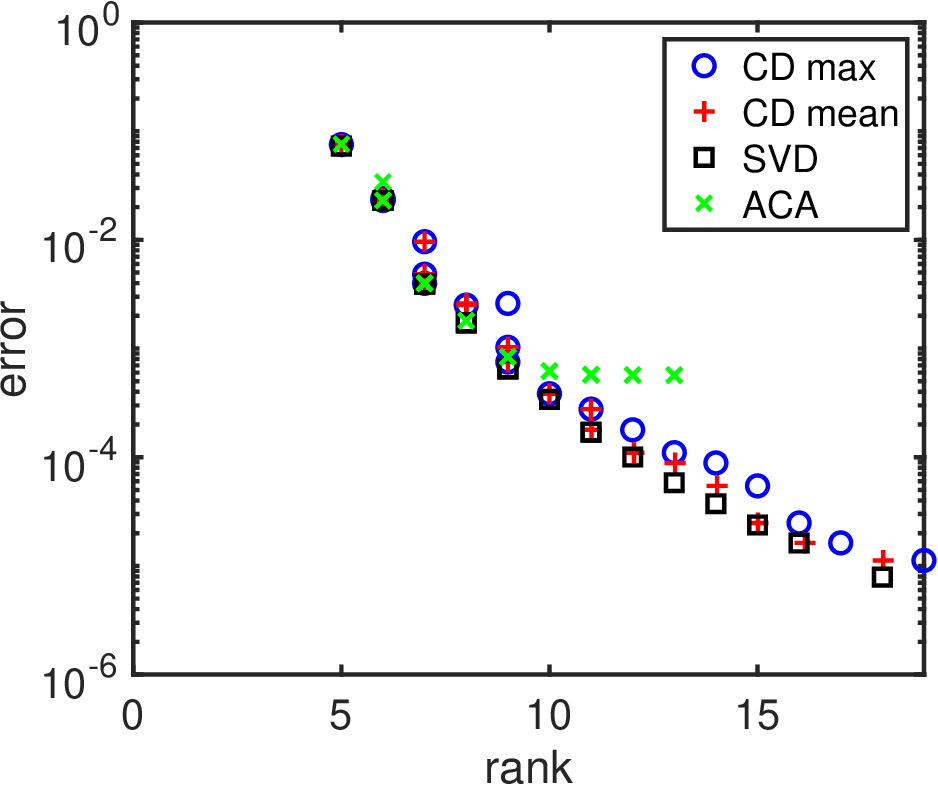}
\includegraphics[width=0.35\textwidth,trim={0.0cm 0.0cm 0.0cm 0.0cm},clip]{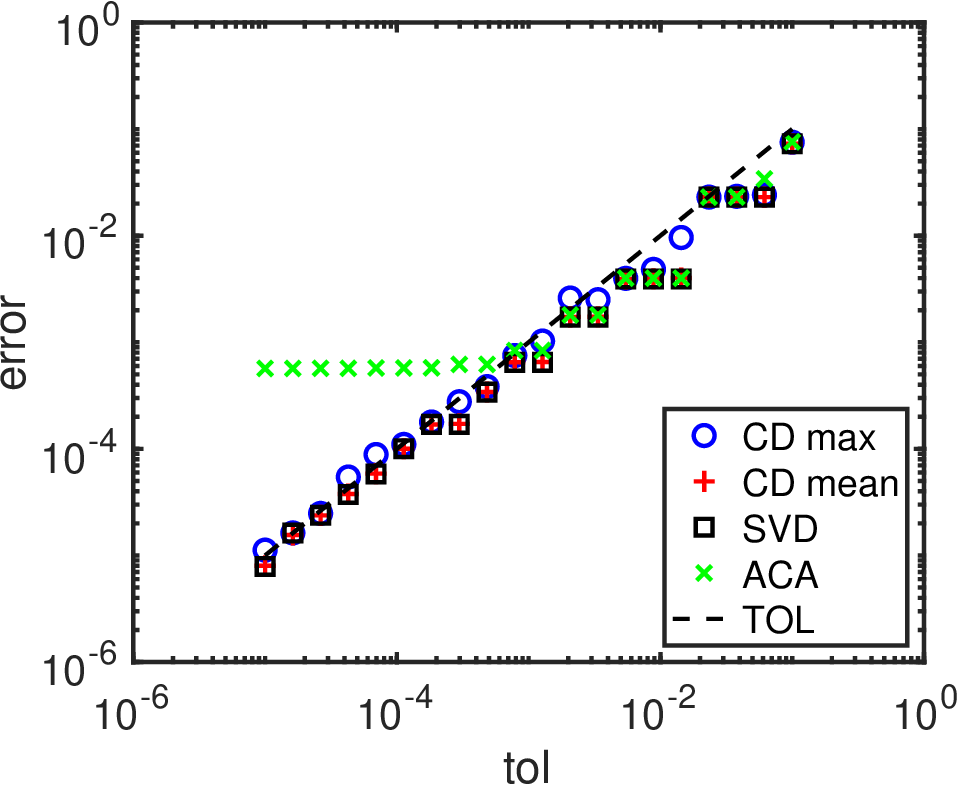}
\includegraphics[width=0.35\textwidth,trim={0.0cm 0.0cm 0.0cm 0.0cm},clip]{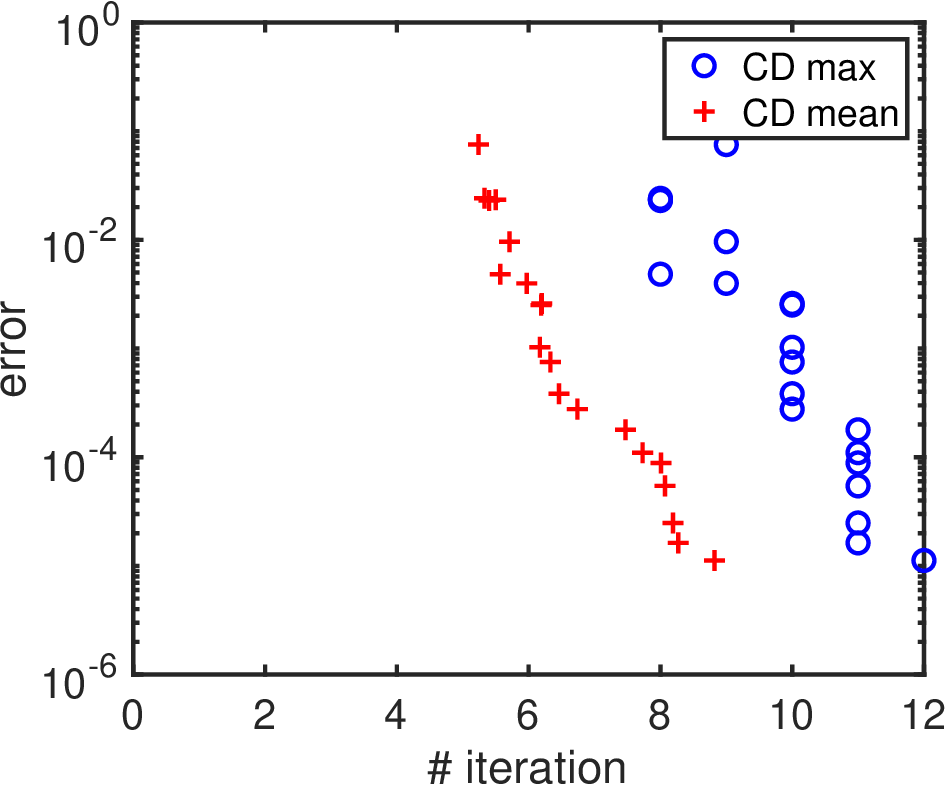}
\includegraphics[width=0.35\textwidth,trim={0.0cm 0.0cm 0.0cm 0.0cm},clip]{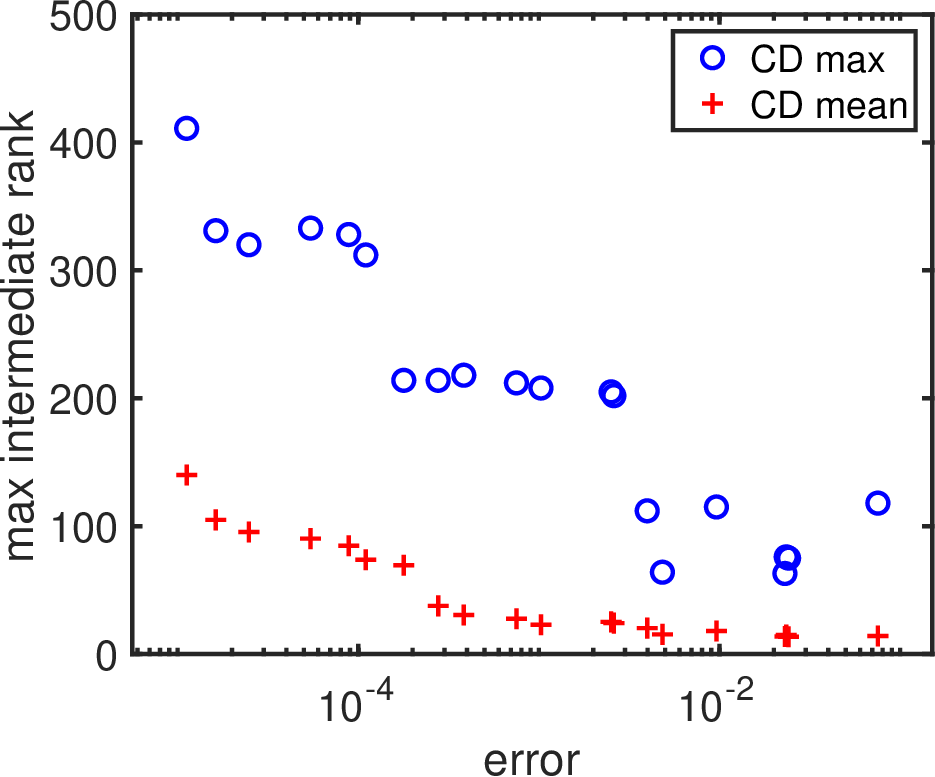}
\caption{Results for approximation of  $G_2$. In the figures, `CD' denotes Cross-DEIM.  \label{fig:G3}}
\end{center}
\end{figure}

The results for $G_2$  with tolerance $\epsilon$ varying from $10^{-5}$ to $10^{-1}$ are reported  
in Figure \ref{fig:G3}. This is a  more challenging example as can be seen from the top two subfigures. The ACA fails to meet the requested tolerance for $\epsilon < O(10^{-3}).$ In particular, ACA deviates from the truncated SVD for ranks above 10. In comparison, Cross-DEIM method gives   results that are close to the truncated SVD for all tolerances and random initial vectors.  From the bottom left figure, we can see the iteration number is not too large for this challenging problem, with the mean ranging from 4-8. Note however that for both the iteration number and max intermediate rank, the max can be much larger that the mean for low tolerances (errors), illustrating some sensitivity of the method to the random initial vector. 

In conclusion, we can see that for both matrices, Cross-DEIM outputs results that are similar to truncated SVD in terms of rank and error. In particular, the numerical error meets the specified tolerance. Compared to an easy problem where the matrix has rapid singular value decay, the more challenging problem requires  a   larger number of iterations and a larger max intermediate rank , which means more computational resources will be needed. It is also clear that for more challenging problem, the computational efficiency of the Cross-DEIM method is more sensitive to the random initial vectors. That's the main motivation to use warm-start in the Cross-DEIM method as will be demonstrated in the next example.

\subsubsection{Parametric matrix approximation}
In this experiment, we consider the approximation of parametric  matrices of size $m \times n$ whose elements are dependent on a parameter $t \in [0,1]$.  Here the discrete time step is $\Delta t = 1/80$. Therefore, we are approximating a sequence of $80$ matrices parametrized by $t$.

We define the matrix elements indexed by $(i,j)$  to be associated with a two-dimensional coordinate 
\[
\begin{pmatrix}
x_i(t) \\ 
y_j(t)
\end{pmatrix} =
\begin{pmatrix}
\phantom{-}\cos(2 \pi t) & \sin (2 \pi t)\\ 
-\sin (2 \pi t) & \cos(2 \pi t) 
\end{pmatrix}\begin{pmatrix}
-1 + ih_x \\ 
-1 + jh_y
\end{pmatrix}, 
\] 
where $h_x = 2/(m+1),  \ \ h_y = 2/(n+1)$. The time dependent coordinates resembles solid body rotation. With this convention, we consider the following two matrices
\begin{gather*}
H_{1}(i,j) = e^{-\left( \left(\frac{x_i}{0.3}\right)^2+\left(\frac{y_j}{0.1}\right)^2\right)},
\ \ \ \ H_2(i,j) = \left(\frac{|x_i+y_j|}{2}\right)^5.
 \end{gather*} 

In both cases, we take $m=n=500,$ and fix the tolerance to be $10^{-2}$. For both cases, we expect the ranks to change according to  the tolerance. We set the remaining Cross-DEIM parameters to be their maximum allowable numbers, effectively inactivating them. 
We consider both the ``cold-" and the ``warm-"start strategies, where the cold-start refers to a random instance of initial vectors $U_0, V_0.$ For the warm-start, we use the singular vectors obtained from the previous time step as initial vectors for approximating $X_n$.

\graphicspath{{figures/parametric_function}}
\begin{figure}[htb]
\begin{center}
\includegraphics[width=0.35 \textwidth,trim={0.0cm 0.0cm 0.0cm 0.0cm},clip]{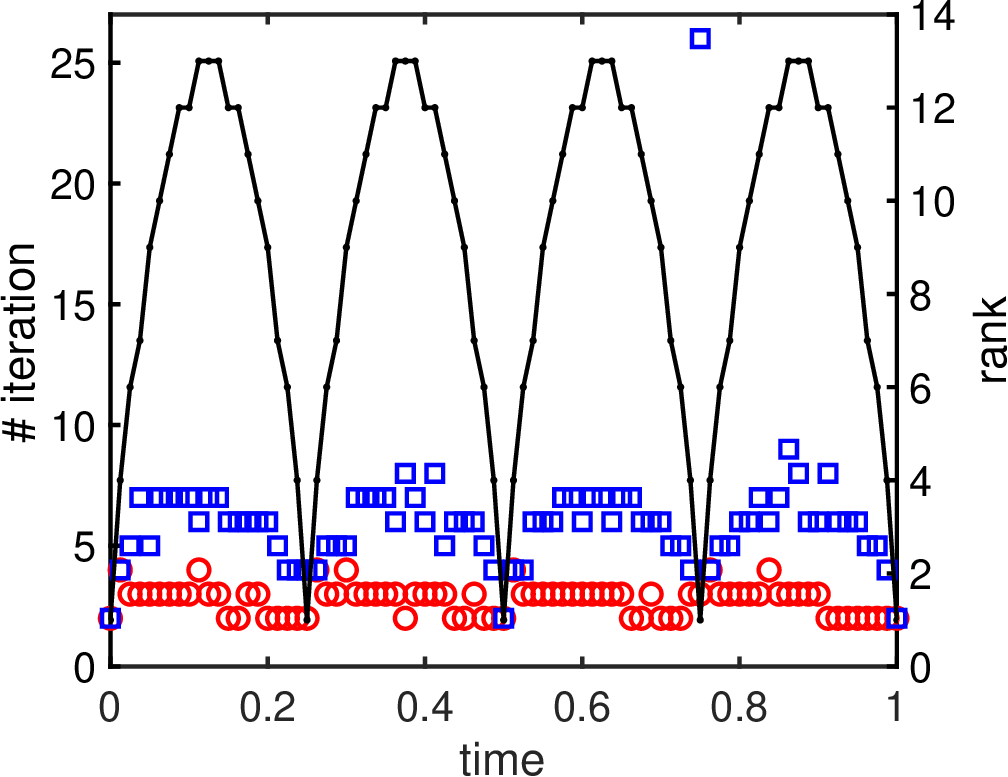}  \ \ \ 
\includegraphics[width=0.35 \textwidth,trim={0.0cm 0.0cm 0.0cm 0.0cm},clip]{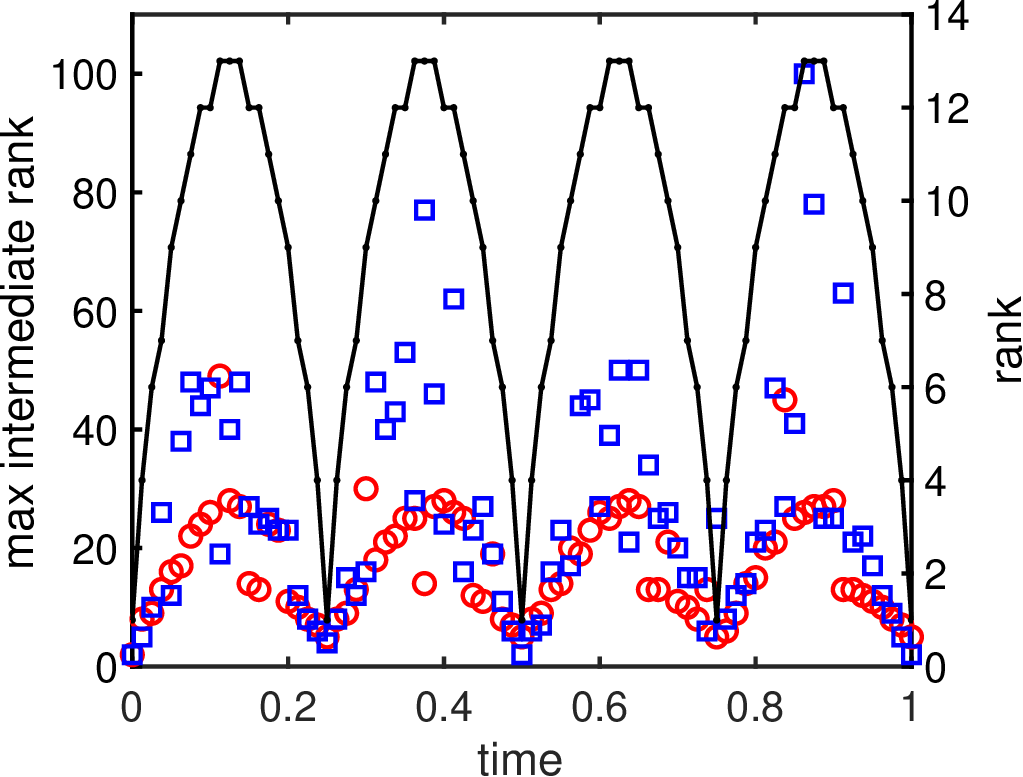}  
\caption{Results for parametric matrix approximation: $H_1$ in terms of iteration number (left figure) and max intermediate rank (right figure). Blue: cold-start, red: warm-start. Black line: numerical rank (the axis is to the right of the figures). \label{fig:parametric1}}
\end{center}
\end{figure}

\begin{figure}[htb]
\begin{center} 
\includegraphics[width=0.35 \textwidth,trim={0.0cm 0.0cm 0.0cm 0.0cm},clip]{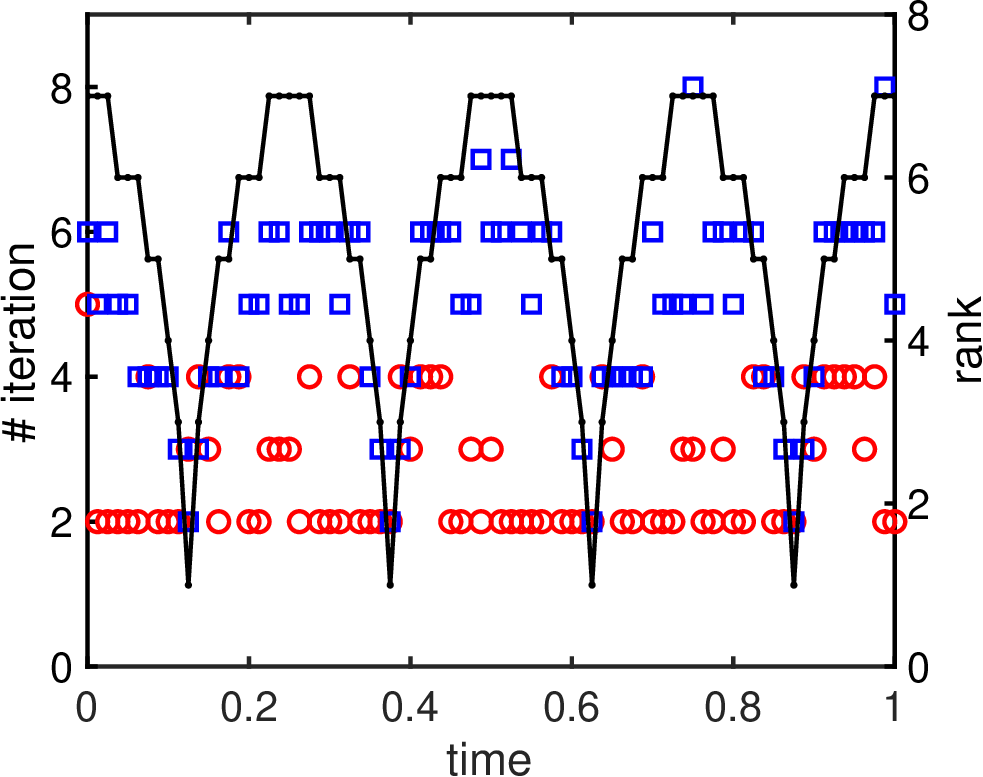} \ \ \ 
\includegraphics[width=0.35 \textwidth,trim={0.0cm 0.0cm 0.0cm 0.0cm},clip]{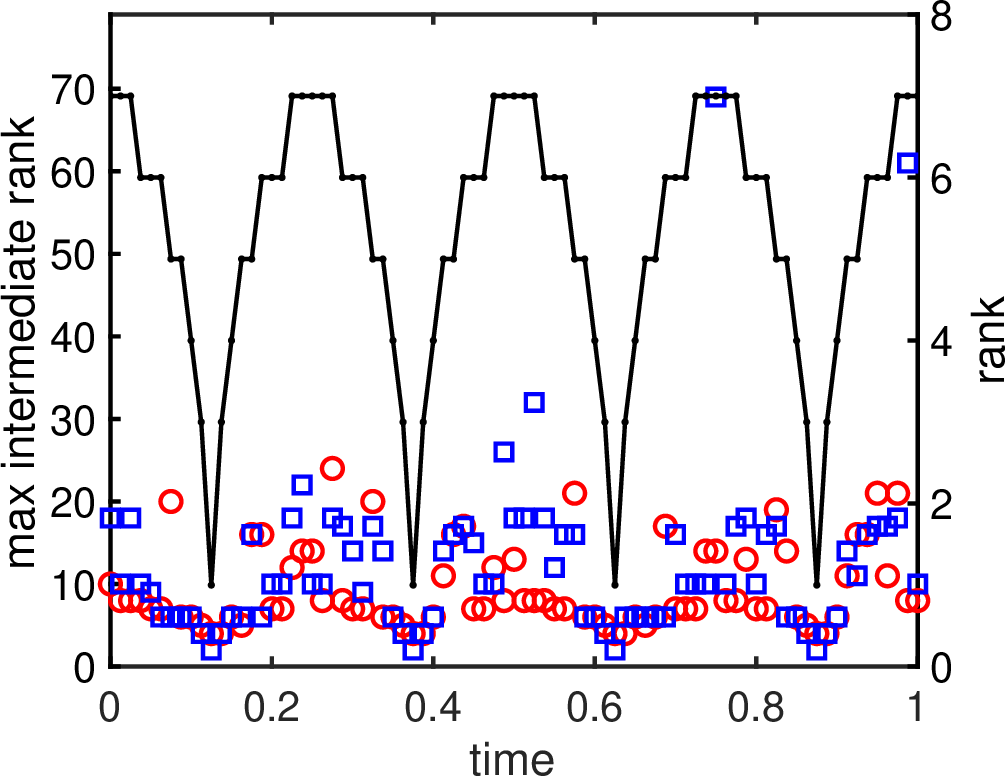}
\caption{Results for parametric matrix approximation: $H_2$ in terms of iteration number (left figure) and max intermediate rank (right figure). Blue: cold-start, red: warm-start. Black line: numerical rank (the axis is to the right of the figures).  \label{fig:parametric2}}
\end{center}
\end{figure}

In Figures \ref{fig:parametric1} and \ref{fig:parametric2}, we compare the the number of iterations for Cross-DEIM to converge (note that the label is to the left in the figure) and the max intermediate rank as a function of the time step. For both functions we see that warm-start strategies work well in bringing down both iteration number of the intermediate ranks (roughly by half). In particular, warm-start removes the sensitivity of computational efficiency on random initial vectors.

\subsection{lrAA examples}
In this subsection, we test the lrAA method. We start with a linear problem from the finite difference approximation to the Laplace's equation to test lrAA without Cross-DEIM. Then we consider nonlinear problems, in particular the  Bratu problem, the fully nonlinear Monge-Amp\'{e}re equations and the time-dependent Allen-Cahn equation. For the nonlinear problems, we   use  lrAA together with Cross-DEIM.

\subsubsection{Laplace's equation}
\graphicspath{{figures/scheduling_laplace}}
We use lrAA to solve Laplace's equation in two dimensions. Since this is a linear problem, we replace Cross-DEIM in Algorithm \ref{algo:AAlr} by the rounding operation in Algorithm \ref{algo:round} to focus on benchmarking the performance of lrAA method by itself.

Precisely, we seek an approximate low-rank solution to   
\[
u_{xx} + u_{yy} = f(x,y), \ \ (x,y) \in [-1,1] \times [-1,1],
\]
with homogenous Dirichlet boundary conditions. To find the approximate solution we discretize this equation using standard second order finite difference approximations for the $x$ and $y$ derivatives. Given an approximation $X(i,j) \approx u(x_i,y_j)$ this results in a function $G_{\Delta}(i,j;X)$ describing the stencil  
\begin{gather}
\begin{split}
G_{\Delta}(i,j;X) = \frac{1}{h_x^2}\left(X(i+1,j) - 2X(i,j) + X(i-1,j)  \right) \\
+ \frac{1}{h_y^2} \left(X(i,j+1) - 2X(i,j) + X(i,j-1)  \right)
\end{split}
\end{gather}
Near the boundaries, some of the terms in this expression will be set to zero to account for the homogenous Dirichlet boundary conditions. Here, note that given the SVD representation of $X$ an element $X(i,j)$ can be obtained at a low arithmetic cost (while maintaining a small memory footprint) by evaluating 
$
X(i,j) = \sum_{k = 1}^{r} \sigma_k u_k(i) v_k(j).
$  

In the numerical examples below, we take the mesh to be  $x_i = -1 + i h_x,$  $h_x = \frac{2}{m+1}$ and $y_j = -1 + j h_y,$ $h_y = \frac{2}{n+1}$. The forcing is chosen as $f(x,y) = -25\exp(-36((x-0.52)^2+(y-0.5)^2))$.

The fixed point function $G(i,j)$ is obtained by applying the pre-conditioned Richardson iteration. We have 
\[
X^{k+1}(i,j) = G(i,j;X^{k},F,\alpha) \equiv  X^k(i,j) + \alpha M(G_{\Delta}(i,j;X^{k}) - F(i,j)),
\] 
where $F(i,j) = f(x_i,y_j)$, $M(\cdot)$ is a preconditioner and $\alpha$ is a parameter. When no preconditioning is used, we choose $\alpha = 0.1 \min(h_x^2,h_y^2)$. We will also use the Exponential Sum (ES) preconditioner based on the approximation to $1/x$ described in \cite{BraHacEXPSUM2005} and when we do, we set $\alpha = 1$. In all experiments we start from a random rank 1 matrix.    

\begin{figure}[htb]
\begin{center}
\includegraphics[width=0.24\textwidth,trim={0.9cm 0.0cm 0.0cm 0.0cm},clip]{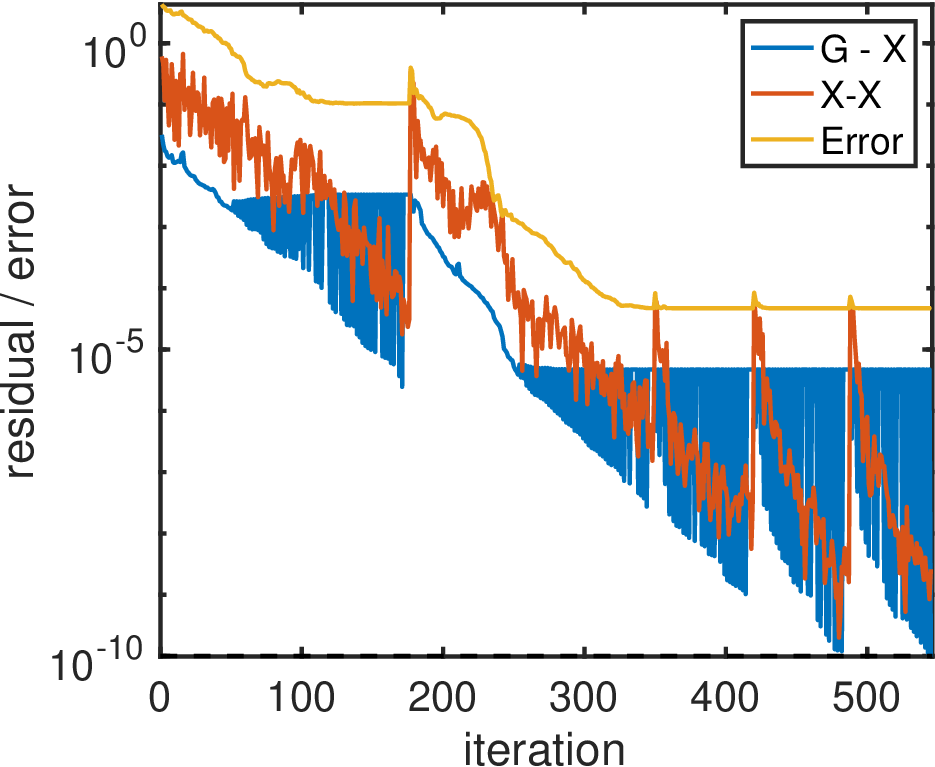}
\includegraphics[width=0.24\textwidth,trim={0.9cm 0.0cm 0.0cm 0.0cm},clip]{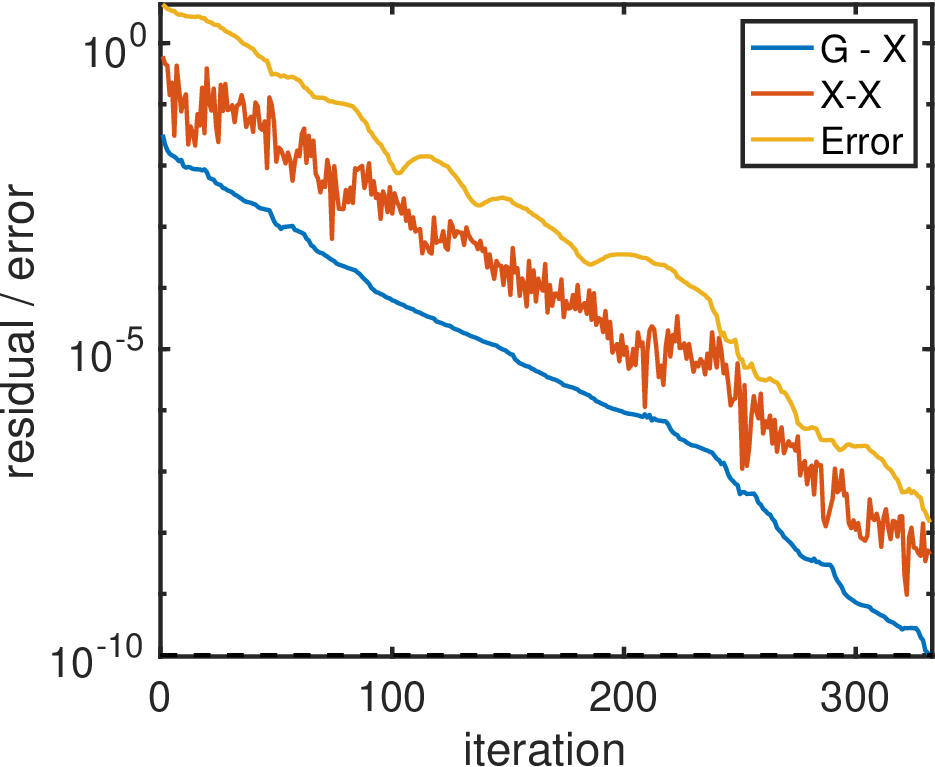}
\includegraphics[width=0.24\textwidth,trim={0.9cm 0.0cm 0.0cm 0.0cm},clip]{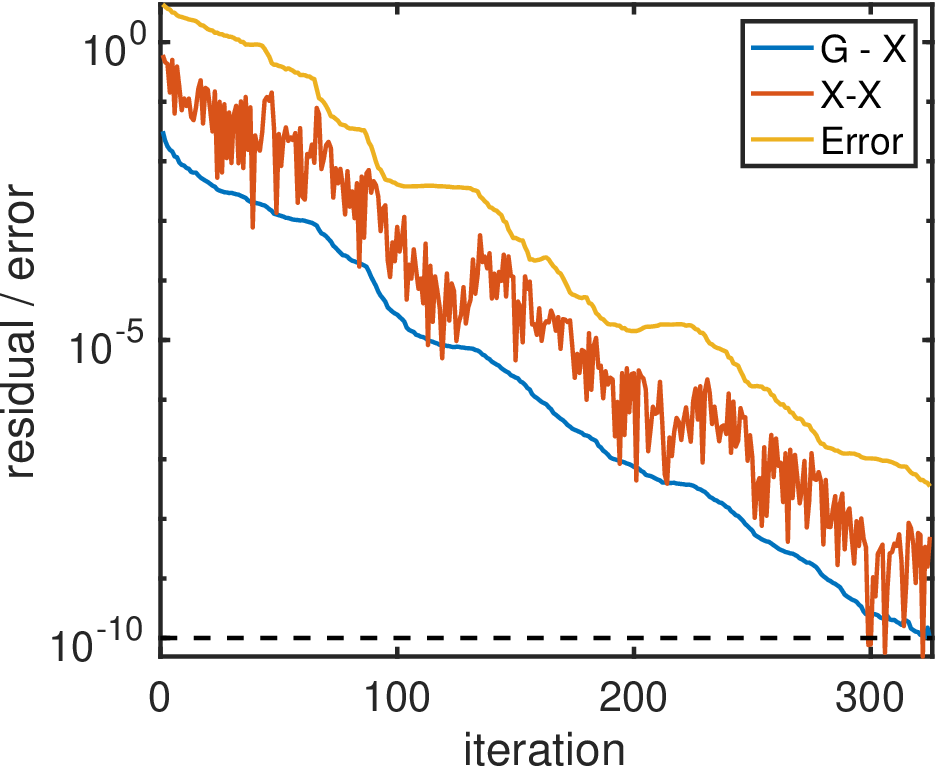}
\includegraphics[width=0.24\textwidth,trim={0.9cm 0.0cm 0.0cm 0.0cm},clip]{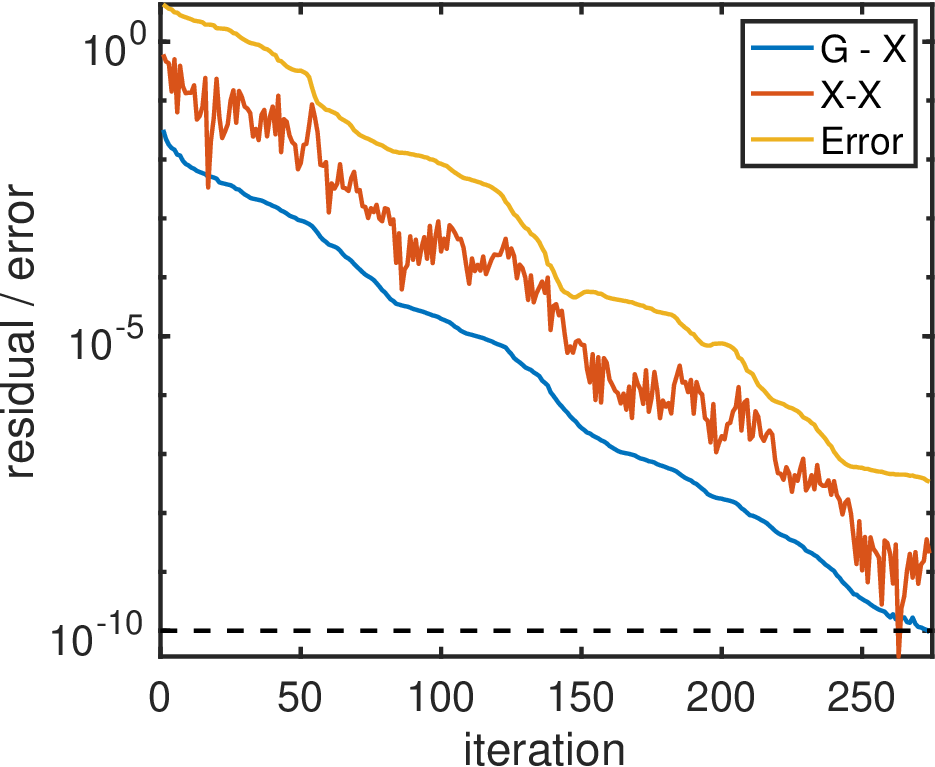}
\includegraphics[width=0.24\textwidth,trim={0.0cm 0.0cm 0.0cm 0.0cm},clip]{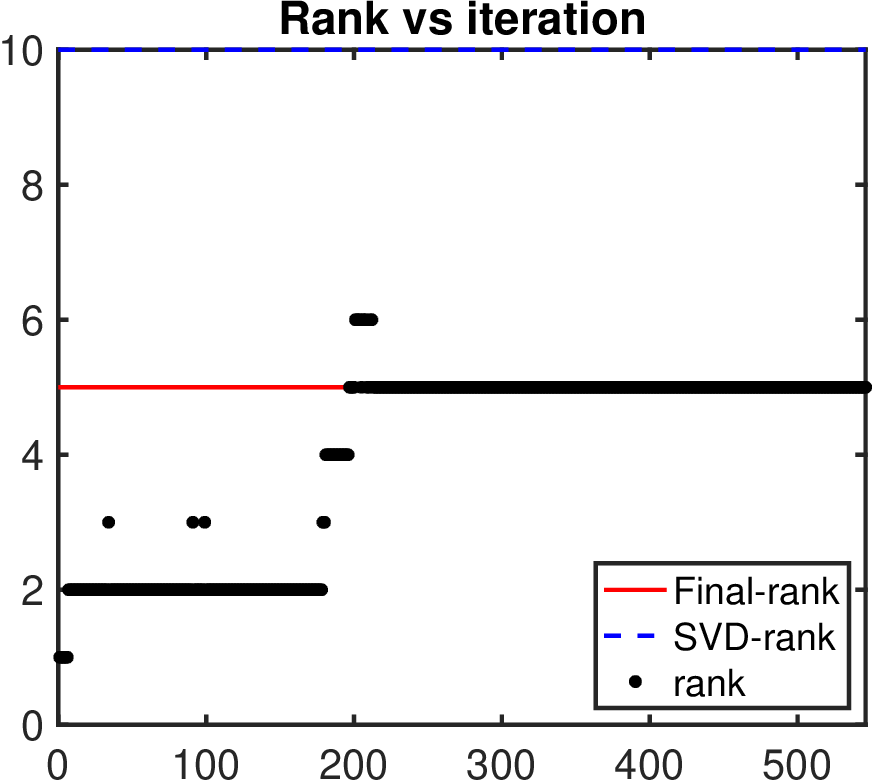}
\includegraphics[width=0.24\textwidth,trim={0.0cm 0.0cm 0.0cm 0.0cm},clip]{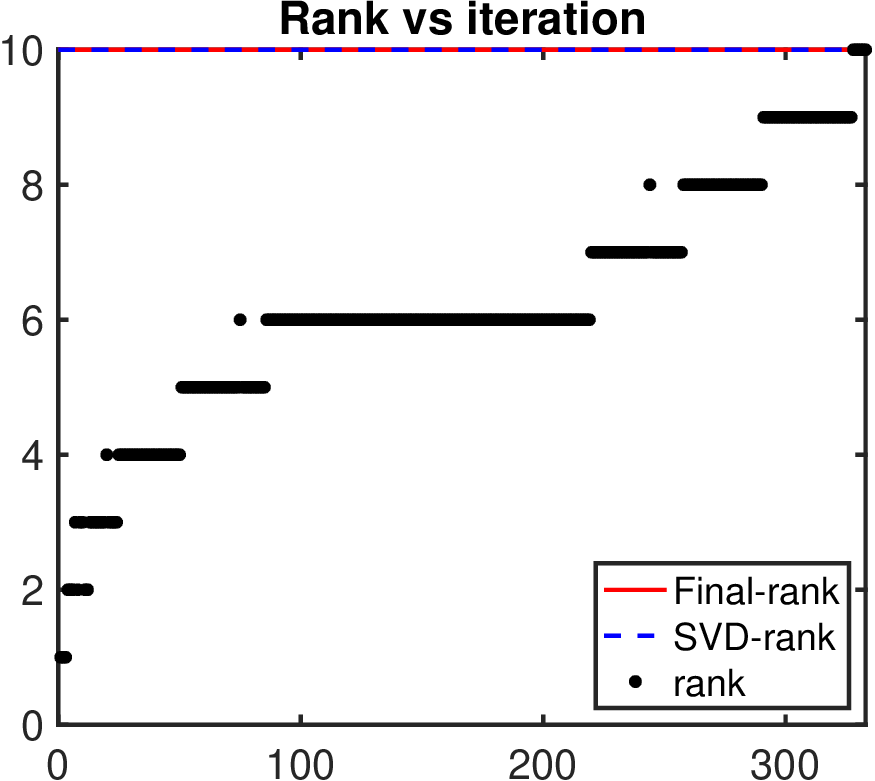}
\includegraphics[width=0.24\textwidth,trim={0.0cm 0.0cm 0.0cm 0.0cm},clip]{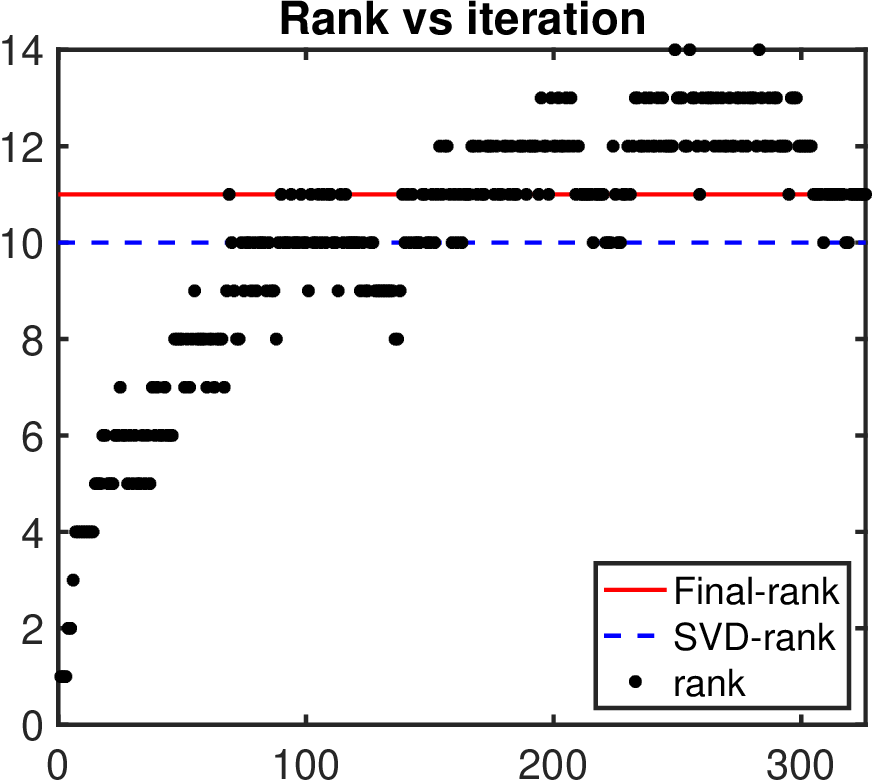}
\includegraphics[width=0.24\textwidth,trim={0.0cm 0.0cm 0.0cm 0.0cm},clip]{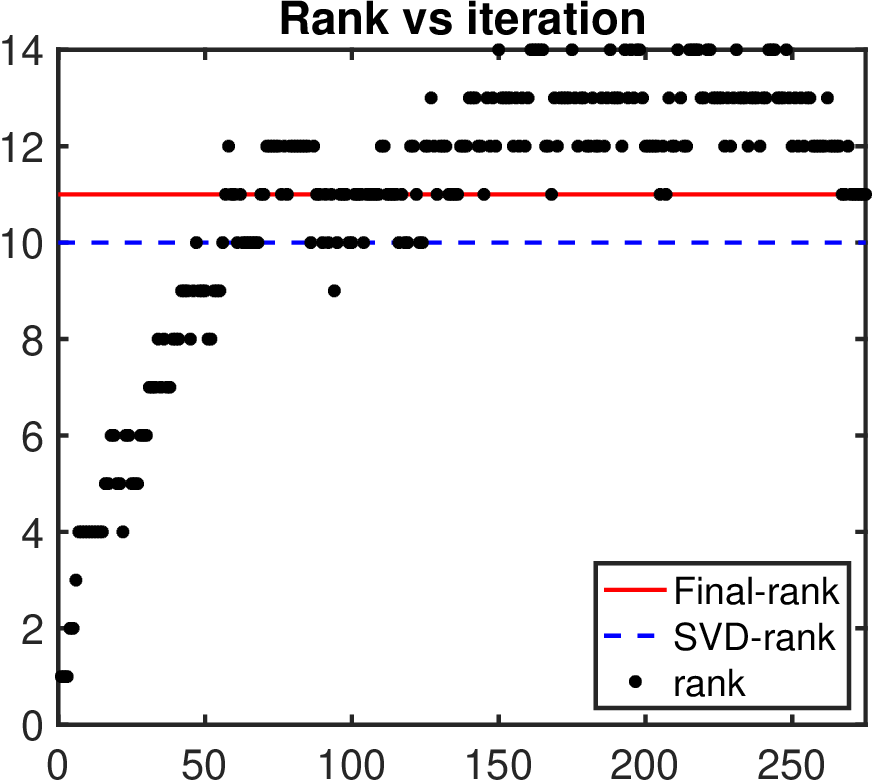}
\caption{Laplace's equation solved by lrAA. In this figure we display the results obtained by the scheduling $\epsilon_{k+1} = \theta \rho_k$. From left to right $\theta$ takes the values 1, 0.5, 0.2, and 0.1. Here all results are obtained with window size size 5 and $m = n = 31$. Results for other grid sizes are qualitatively similar. \label{fig:schedule1}}
\end{center}
\end{figure}

\begin{figure}[htb]
\begin{center}
\includegraphics[width=0.235\textwidth,trim={0.0cm 0.0cm 0.0cm 0.0cm},clip]{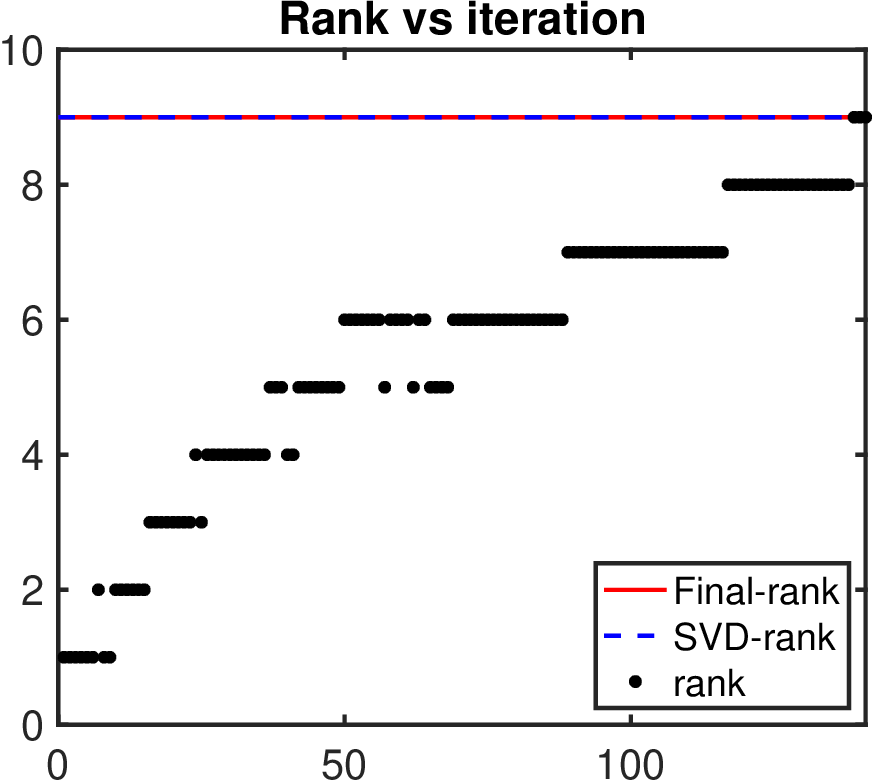}
\includegraphics[width=0.235\textwidth,trim={0.0cm 0.0cm 0.0cm 0.0cm},clip]{scl_rank_5_5_2}
\includegraphics[width=0.235\textwidth,trim={0.0cm 0.0cm 0.0cm 0.0cm},clip]{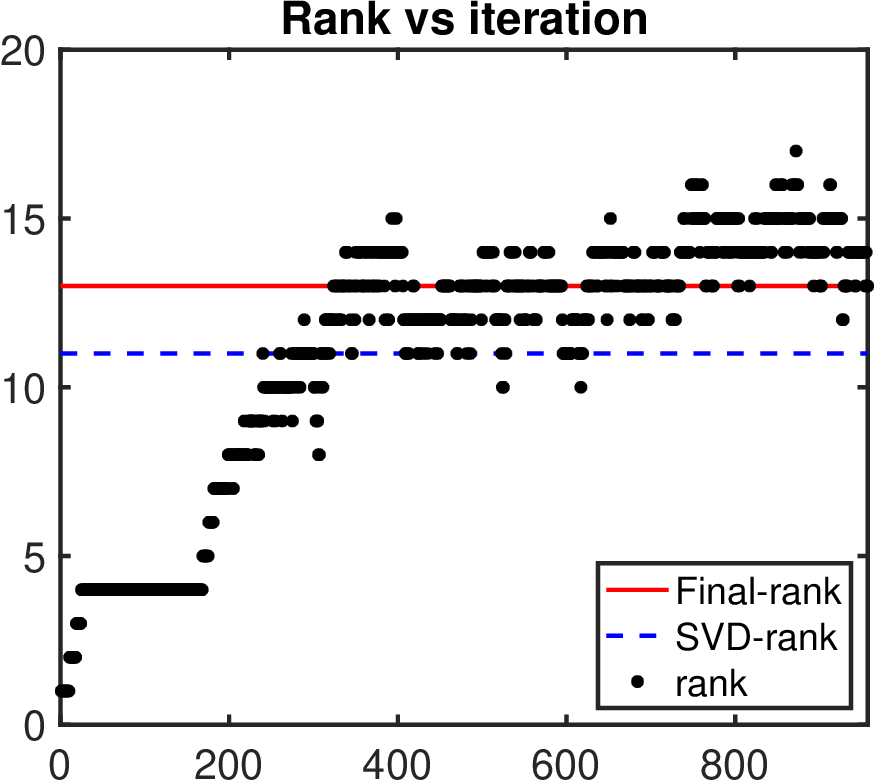}
\includegraphics[width=0.265\textwidth,trim={0.5cm 0.2cm 0.0cm 0.0cm},clip]{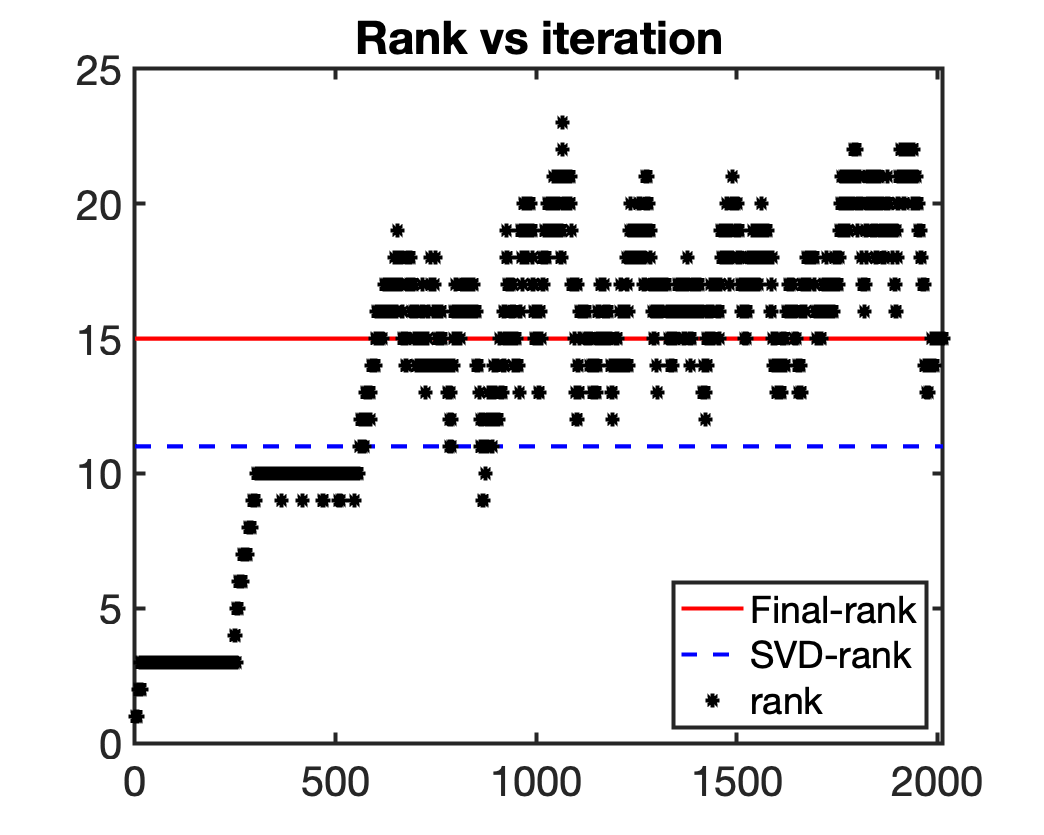}
\caption{Laplace's equation solved by lrAA. In this figure we display the intermediate ranks using the scheduling $\epsilon_{k+1} = 0.5 \rho_k$. From left  to right we have $m=n = 15, 31, 63, 127$. \label{fig:schedule2}}
\end{center}
\end{figure}

\subsubsection*{No preconditioner} We first take $M=I,$ i.e. no preconditioner is used in the Richardson iteration.
In Figure \ref{fig:schedule1}, we study the effect of scheduling of the rounding tolerance for $X$ and $G$. We   consider the impact on the rate of convergence and the intermediate rank of the solution  of different choices of $\theta$. We take the window size $\hat{m}=5$ and choose  $n=m = 31$. 
In Figure \ref{fig:schedule1}, we display the error, $\rho_k$ (denoted {\sf G-X})  and $ \| X_{k+1} - X_{k} \|$ (denoted {\sf X-X}) as  functions of the iteration count. We also display the rank during the iteration and the final rank along with the rank obtained by truncating the SVD of the solution computed using a direct method (Gauss elimination) at the stopping tolerance ${\rm TOL}=10^{-10}$ (denoted by SVD-rank). From left to right in Figure \ref{fig:schedule1}, $\theta$ takes the values 1, 0.5, 0.2, and 0.1. The choice $\theta = 1$ makes the residual and rounding equal, and as can be seen in the figure this does not lead to a robust scheduling. Here we find that $0.1 \le \theta \le 0.5$ appears to good convergence and also keeps the rank of the solution iterates low and close to the rank of the converged solution. In particular, at this mesh size, the more aggressive truncation with $\theta=0.5$ yields a monotonically increasing numerical rank till convergence is reached, while the smaller $\theta=0.2, 0.1$ give slight overshoot of the numerical rank during iteration. We also see that $ \| X_{k+1} - X_{k} \|$ is much more  oscillatory  than   $\rho_k$ and we therefore always use $\rho_k$ for scheduling the rounding tolerance and for checking the stopping criteria.  In Figure \ref{fig:schedule2} we use $\theta = 0.5$ and report the intermediate numerical rank as a function of iteration count for $m=n = 15, 31, 63, 127$. The iteration is terminated when the residual is smaller than $10^{-10}$. As can be seen, the scheduling keeps the intermediate ranks in control. However, for  finer mesh size,  the intermediate rank has some overshoot. Nevertheless, the small but apparent effect of grid size indicates a more advanced adaptive scheduling could give yet better   results.

\begin{figure}[htb]
\begin{center}
\includegraphics[width=0.325\textwidth,trim={0.cm 0.0cm 0.9cm 0.0cm},clip]{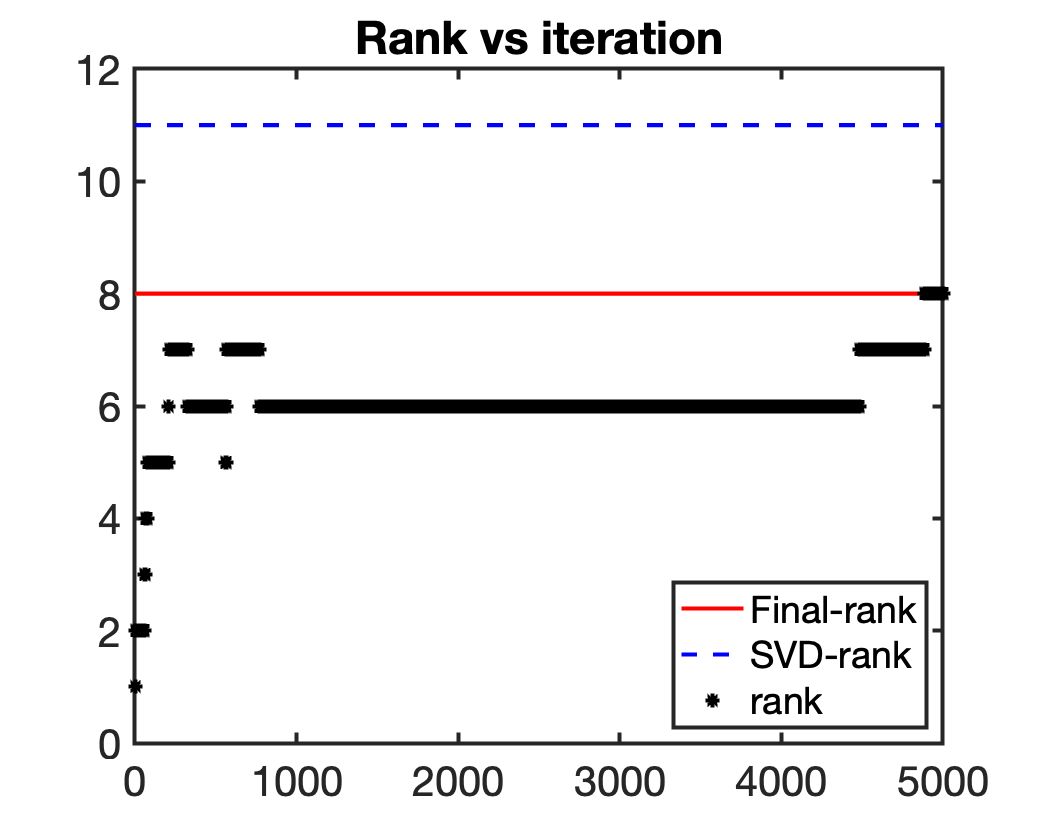}\ \ \ 
\includegraphics[width=0.30\textwidth,trim={0.cm 0.0cm 0.0cm 0.0cm},clip]{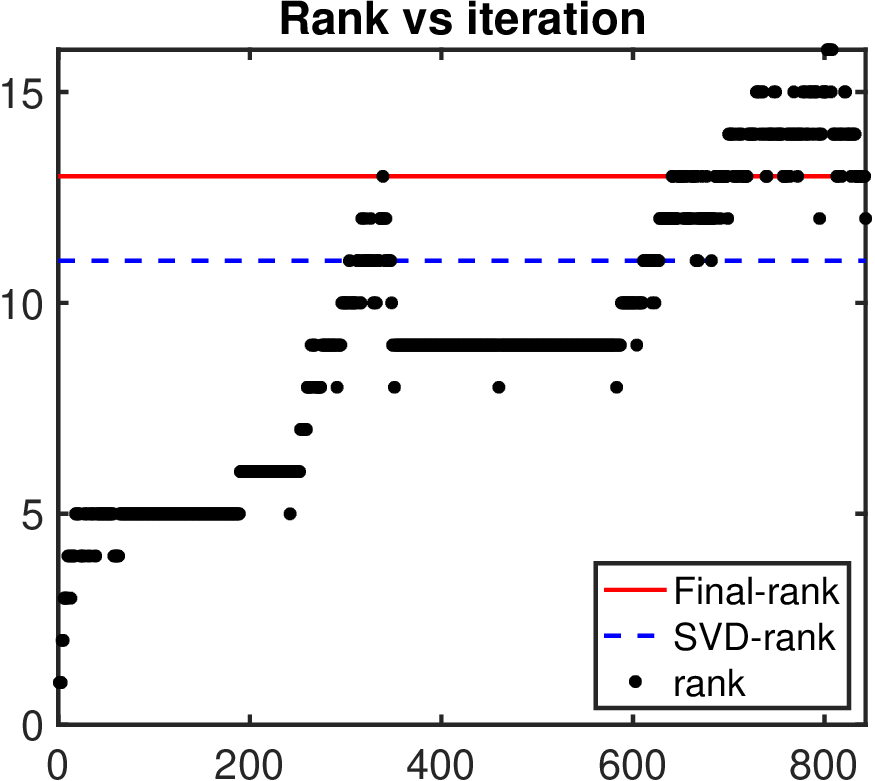}\ \ \ 
\includegraphics[width=0.30\textwidth,trim={0.cm 0.0cm 0.0cm 0.0cm},clip]{scl_rank_6_5_2}
\includegraphics[width=0.32\textwidth,trim={0.cm 0.0cm 1.35cm 0.0cm},clip]{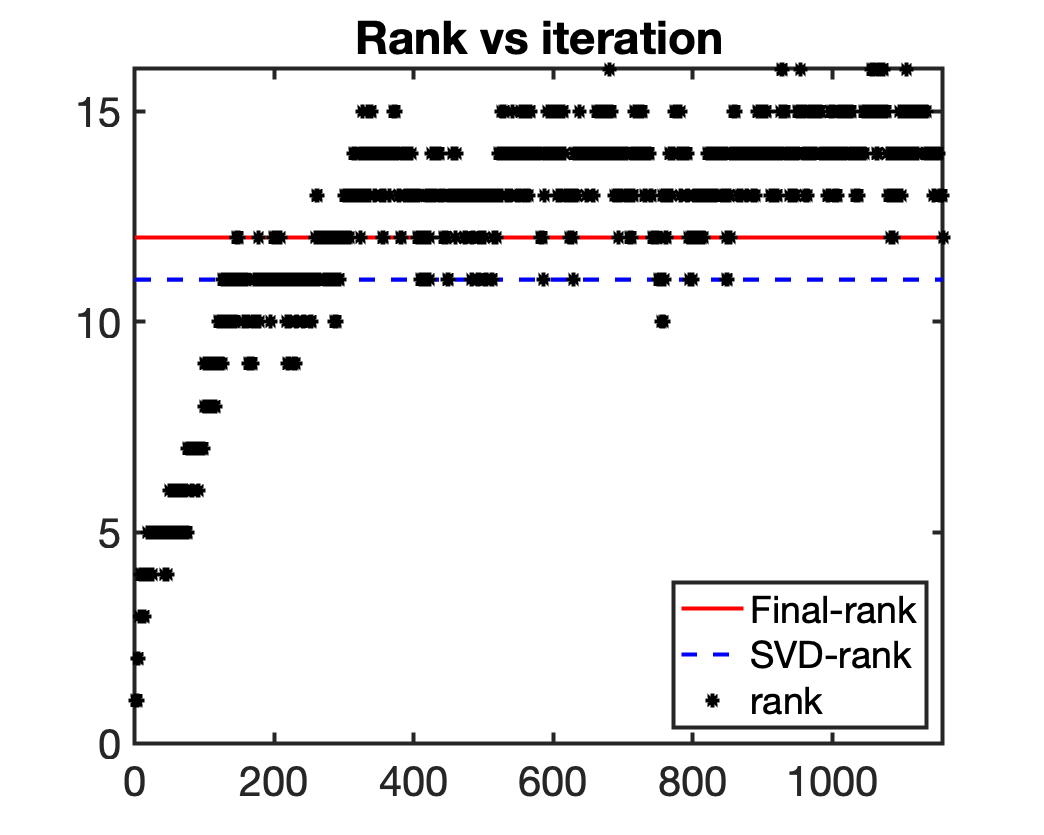}
\includegraphics[width=0.32\textwidth,trim={0.cm 0.0cm 1.35cm 0.0cm},clip]{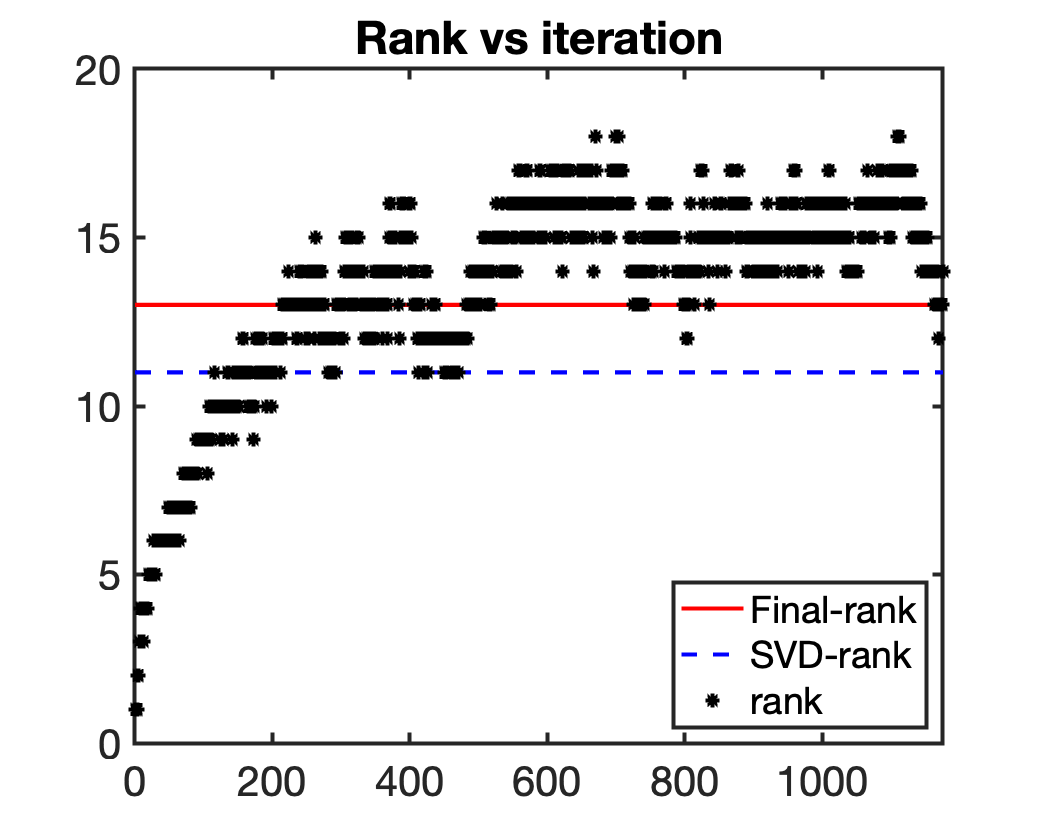}
\includegraphics[width=0.32\textwidth,trim={0.cm 0.0cm 1.35cm 0.0cm},clip]{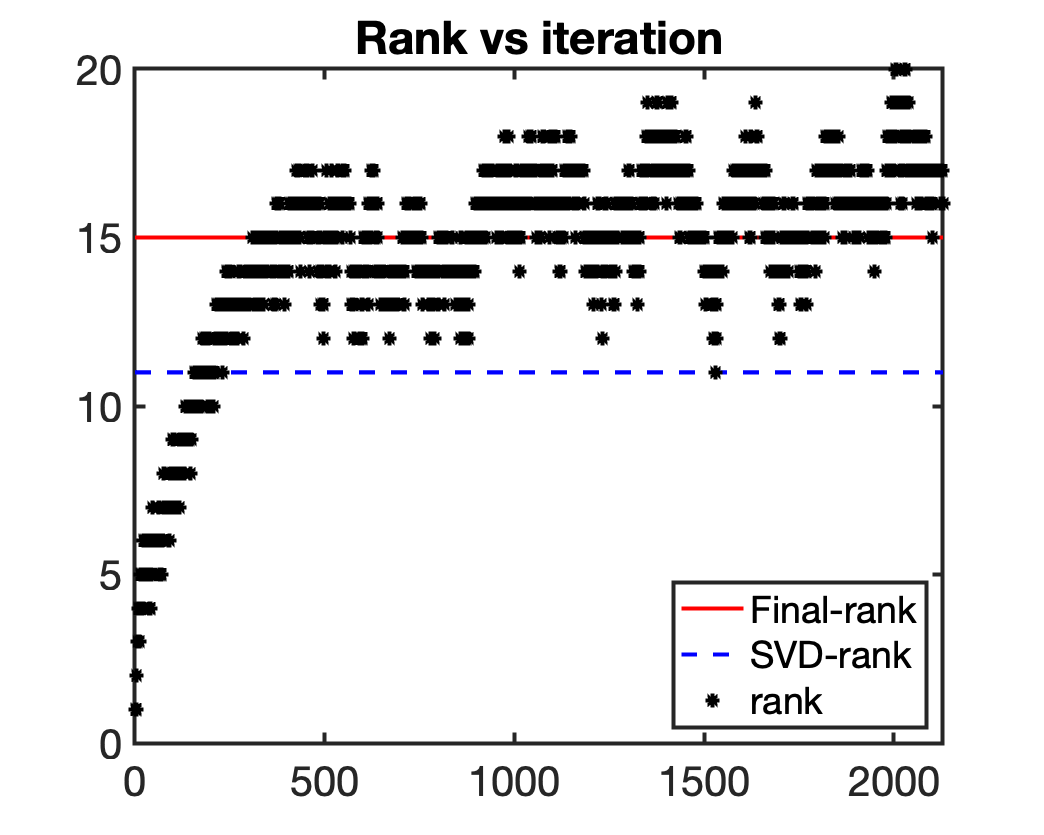}
\caption{Laplace's equation solved by lrAA. This figure displays the intermediate rank as a function of iteration count for window size 1,3,5,7,10,20 (top left to bottom right). We use the scheduling $\epsilon_{k+1} = 0.5 \rho_k$, and $n=m=63$. \label{fig:schedule5}}
\end{center}
\end{figure}

In Figure \ref{fig:schedule5}, we study the effect of different window size $\hat{m}.$ We display the intermediate rank as a function of iteration count for window sizes 1, 3, 5, 7, 10, 20. The iteration is stopped when $\rho_k$ is below ${\rm TOL}=10^{-10}$ or when the number of iterations have reached 5000. This computation is done with $m=n=63$. It is clear from the results that lrAA requires larger number of iterations when the window size is too small or too large. We have found that a window size equal to 5 works well and in all the examples below we use this window size unless explicitly otherwise noted.

\begin{figure}[htb]
\begin{center}
\includegraphics[width=0.24\textwidth,trim={0.93cm 0.85cm 0.0cm 0.0cm},clip]{scl_res_5_5_2}
\includegraphics[width=0.23\textwidth,trim={0.0cm 0.0cm 0.0cm 0.0cm},clip]{scl_rank_5_5_2}
\includegraphics[width=0.24\textwidth,trim={0.93cm 0.85cm 0.0cm 0.0cm},clip]{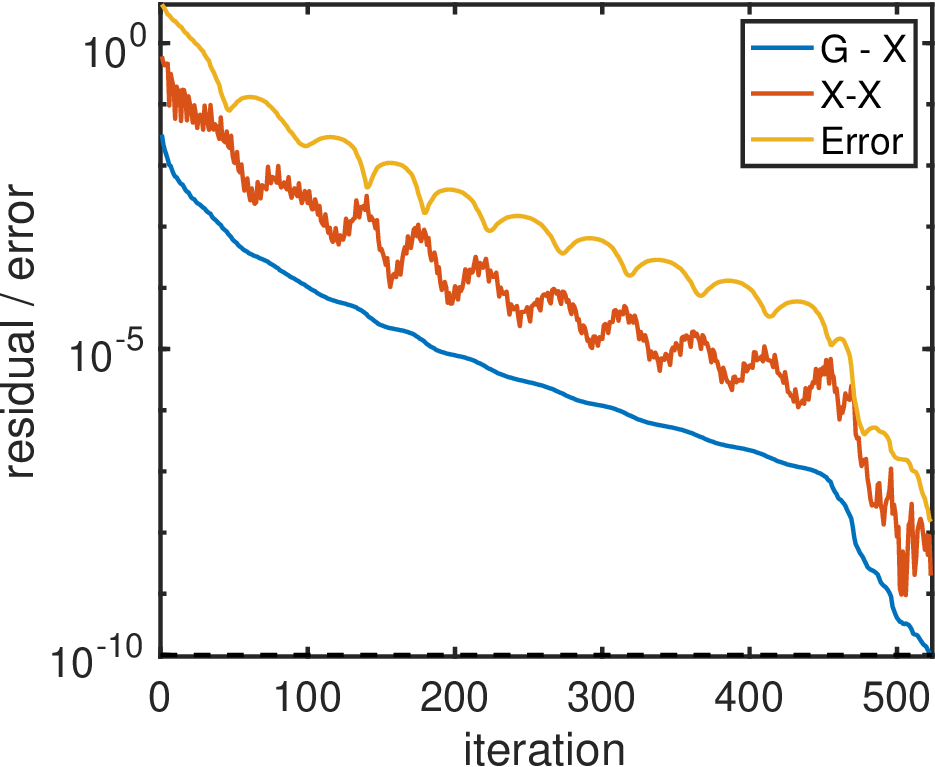}
\includegraphics[width=0.23\textwidth,trim={0.0cm 0.0cm 0.0cm 0.0cm},clip]{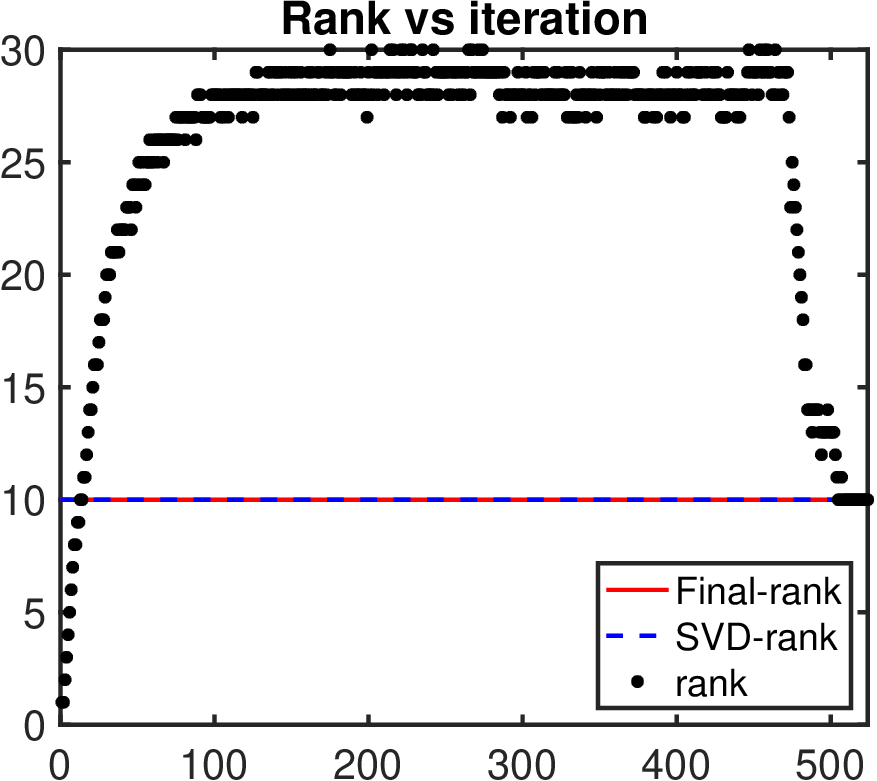}
\caption{Laplace's equation solved by lrAA. The two left graphs show the residuals, errors (left) and ranks (right) obtained using scheduling $\epsilon_{k+1} = 0.5 \rho_k$. The two right most graphs show the same but with a fixed truncation set at the level of the iteration stopping tolerance,  $\epsilon_{k+1} = 10^{-10}$. Here all results are for $n=m=31$. \label{fig:schedule3}}
\end{center}
\end{figure}

  In Figure \ref{fig:schedule3}, we demonstrate the importance of scheduling, i.e. the need to vary the truncation tolerance $\epsilon_{k}$ with the iteration $k$. We display the residuals and ranks for computations with $n = m = 32$   with scheduling $\epsilon_{k+1} = 0.5 \rho_k$, and without, i.e. using a fixed tolerance $\epsilon_{k+1} = 10^{-10}, \forall k$. It is clear that without scheduling, the numerical rank reaches full rank early in the iteration (at about iteration 10, the numerical rank already surpassed the final rank). Therefore, without scheduling, the low-rank methods are pointless.  These computations illustrate the essential importance of a scheduled rounding to keep the numerical  ranks small throughout the iteration. 
  
\begin{figure}[htb]
\begin{center}
\includegraphics[width=0.24\textwidth,trim={0.93cm 0.0cm 0.0cm 0.0cm},clip]{scl_res_5_5_2}
\includegraphics[width=0.24\textwidth,trim={0.93cm 0.0cm 0.0cm 0.0cm},clip]{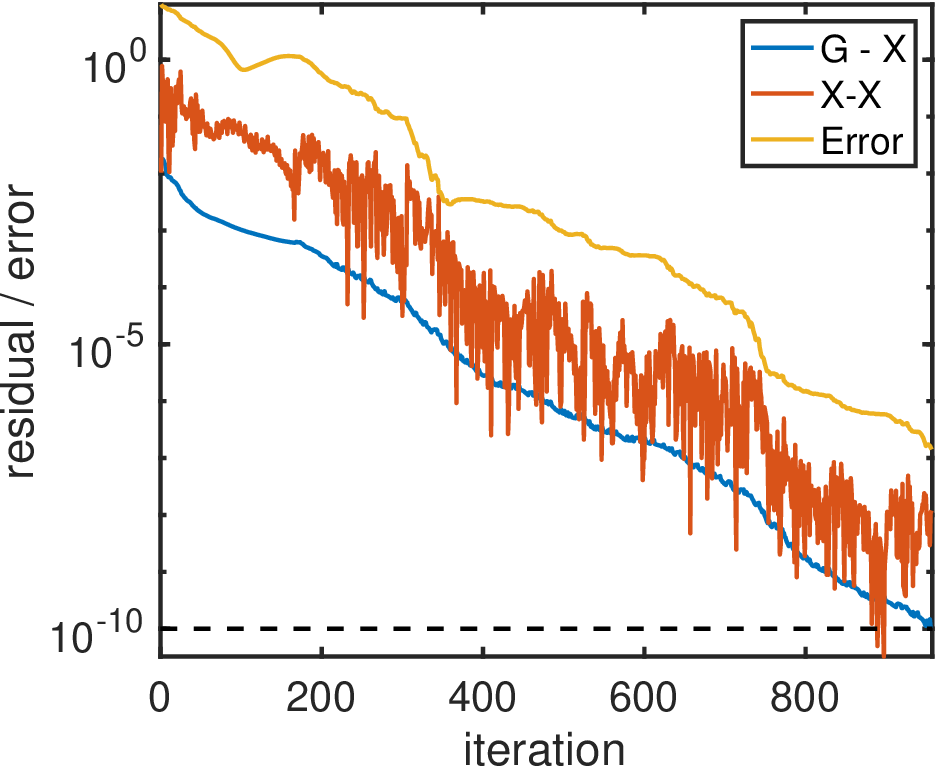}
\includegraphics[width=0.24\textwidth,trim={0.93cm 0.0cm 0.0cm 0.0cm},clip]{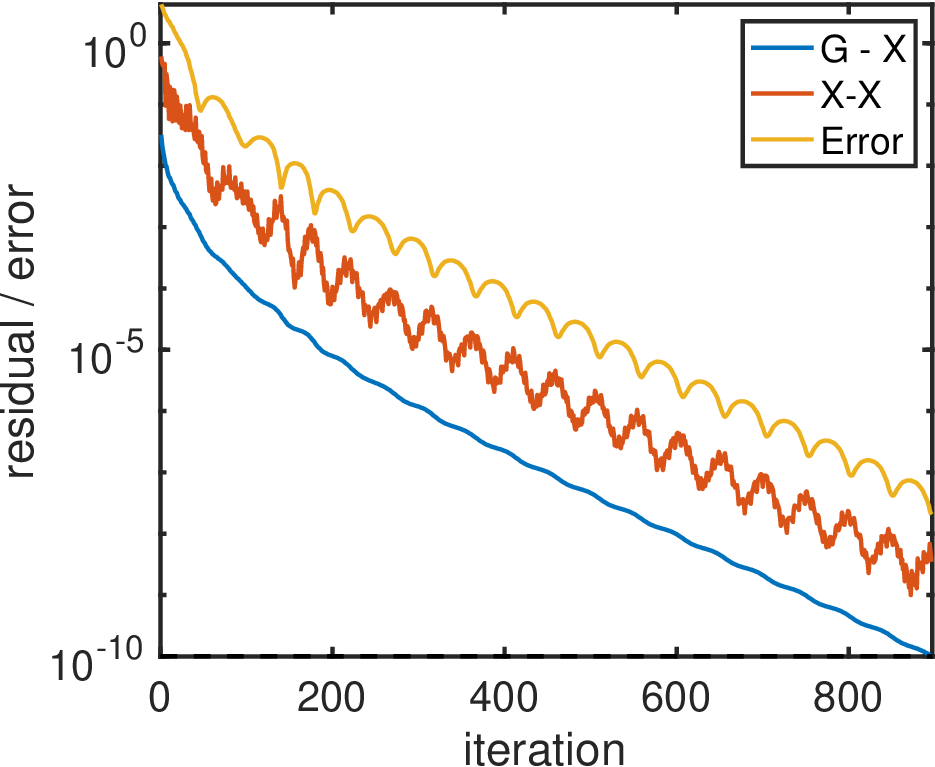}
\includegraphics[width=0.24\textwidth,trim={0.93cm 0.0cm 0.0cm 0.0cm},clip]{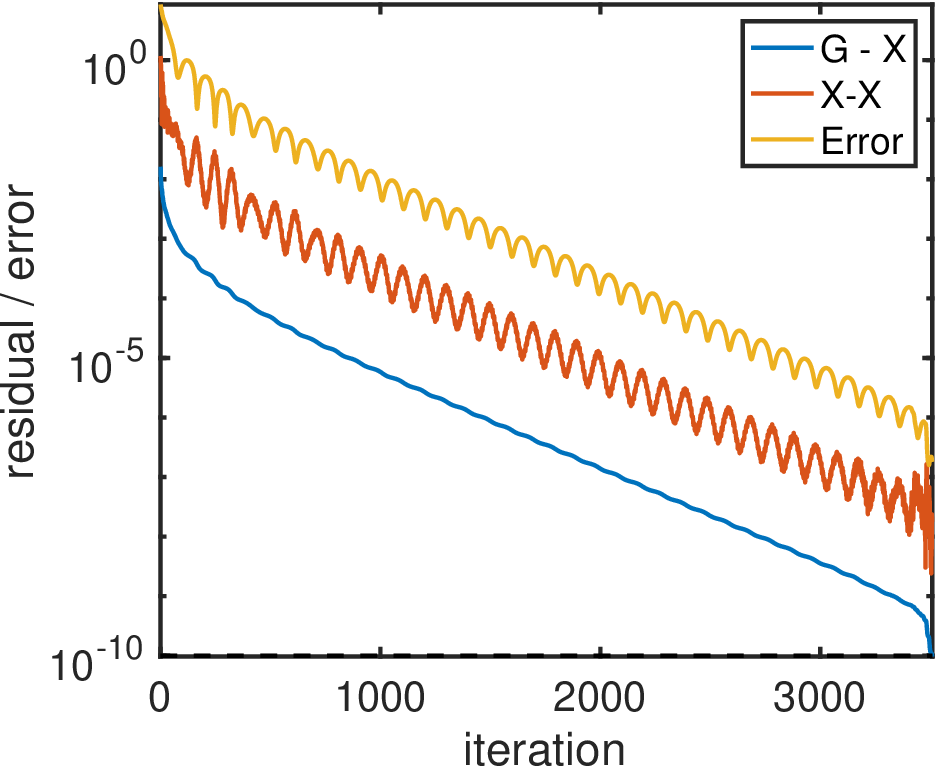}
\caption{Laplace's equation solved by lrAA and full-rank AA. The two leftmost graphs show the residuals and errors obtained by lrAA using scheduling $\epsilon_{k+1} = 0.5 \rho_k$, for $n=m=31$ (left) and 63 (right). The two rightmost graphs show the residuals and errors obtained using the standard full-rank AA, for $n=m=31$ (left) and 63 (right).  \label{fig:schedule4}}
\end{center}
\end{figure}

  In Figure \ref{fig:schedule4}, we demonstrate a comparison of lrAA with the full-rank AA in terms of iteration number. Here, all parameters ($m, n, {\rm TOL}, \hat{m}$) are the same between the two methods, and we consider two different mesh size $n=m=31$ and $n=m=63.$ By comparison, we can see that lrAA has a smaller iteration number compared to full-rank AA. For  $n=m=31$, lrAA's iteration number is about one half of full-rank AA. For the finer mesh, $n=m=63$, lrAA's iteration number is about one quarter of full-rank AA.

\subsubsection*{Preconditioning by Exponential Sums} In this part, we use a well-known low-rank preconditioner: Exponential Sum (ES), to enable calculations for fine mesh size.
When preconditioning low-rank methods, the preconditioner is only allowed to operate on the low rank factors of the solution \cite{grasedyck2013literature,bachmayr2023low}. 
In the ES preconditioner, this is achieved by approximating the inverse as a sum of Kronecker products, \cite{hackbusch2006low,hackbusch2006low2}. 

For problems where the Laplacian is the principal part of the PDE, the ES preconditioner can be very effective.   If the SVD of the residual is  $U_{R}S_{R}V_{R}^T = G_{\Delta}(X^{k}) - F,$ then the residual is applied as follows 
\[
M(G_{\Delta}(X^{k}) - F) = M(U_{R}S_{R}V_{R}^T) = - \sum_{k=1}^{n_{\rm ES}} \alpha_k (e^{\beta_k D_{xx}} U_{R}) S_{R} (e^{\beta_k D_{yy}} V_{R})^T.   
\]
Here $D_{xx}$ and $D_{yy}$ are the one dimensional finite difference matrices for approximating the second derivative in $x$ and $y$. The parameters $\alpha_k$ and $\beta_k$ are chosen so that the preconditioner is a good approximate inverse. Here we exclusively use the weights from \cite{BraHacEXPSUM2005} available from the repository \url{https://gitlab.mis.mpg.de/scicomp/EXP_SUM}. Note that the application of ES preconditioner requires sum of low-rank matrices. This can be computed by rounding or Cross-DEIM as well. To save space, we do not report the details in this paper with regard to this operation.

\begin{figure}[htb]
\begin{center}
\includegraphics[width=0.24\textwidth,trim={0.9cm 0.75cm 0.0cm 0.0cm},clip]{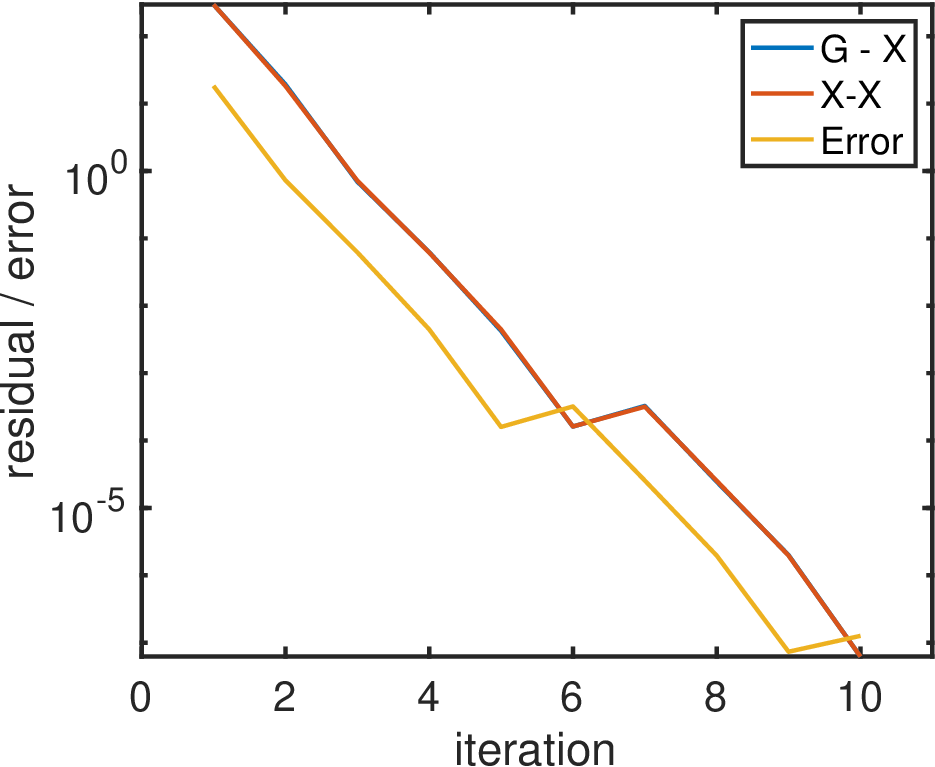}  
\includegraphics[width=0.23\textwidth,trim={0.cm 0.0cm 0.0cm 0.0cm},clip]{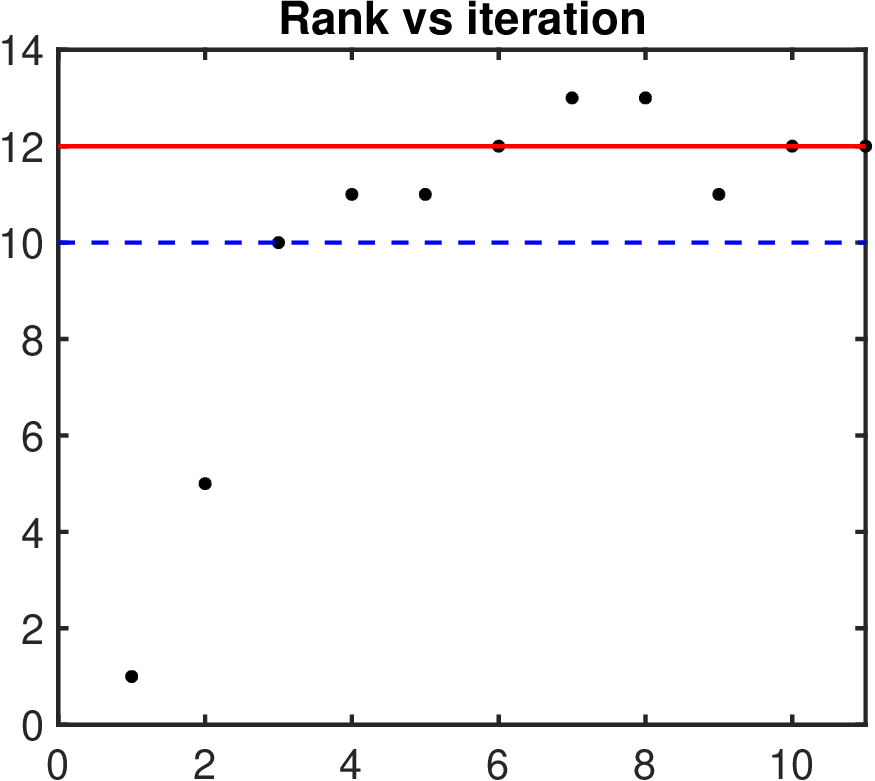}  
\includegraphics[width=0.24\textwidth,trim={0.9cm 0.75cm 0.0cm 0.0cm},clip]{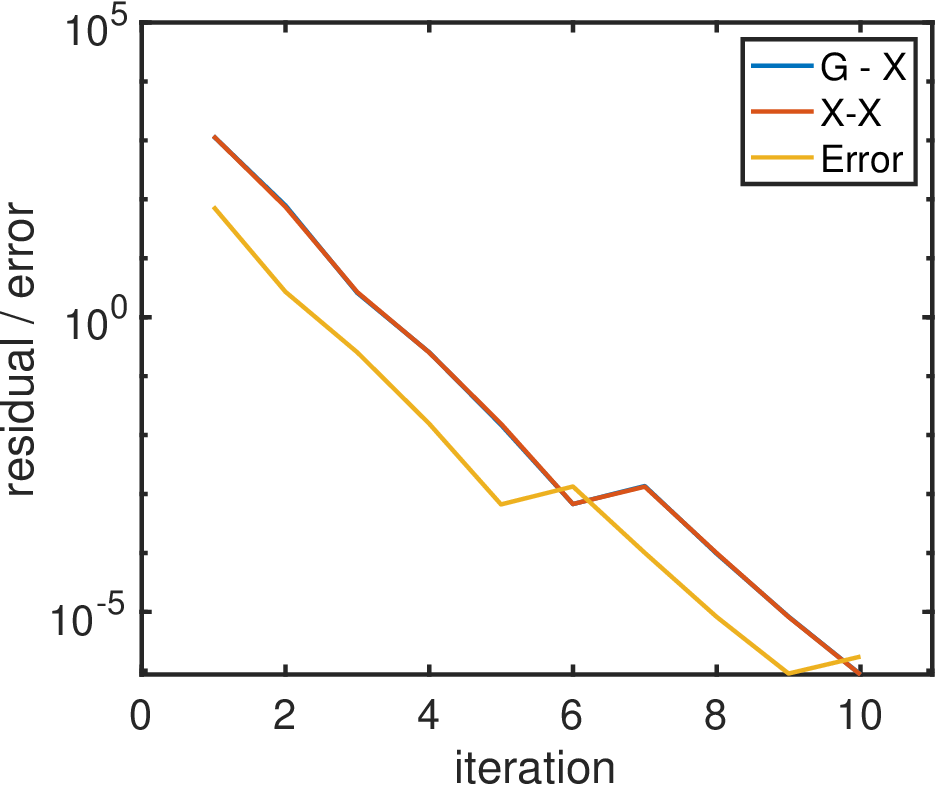} \ \ 
\includegraphics[width=0.225\textwidth,trim={0.cm 0.0cm 0.0cm 0.0cm},clip]{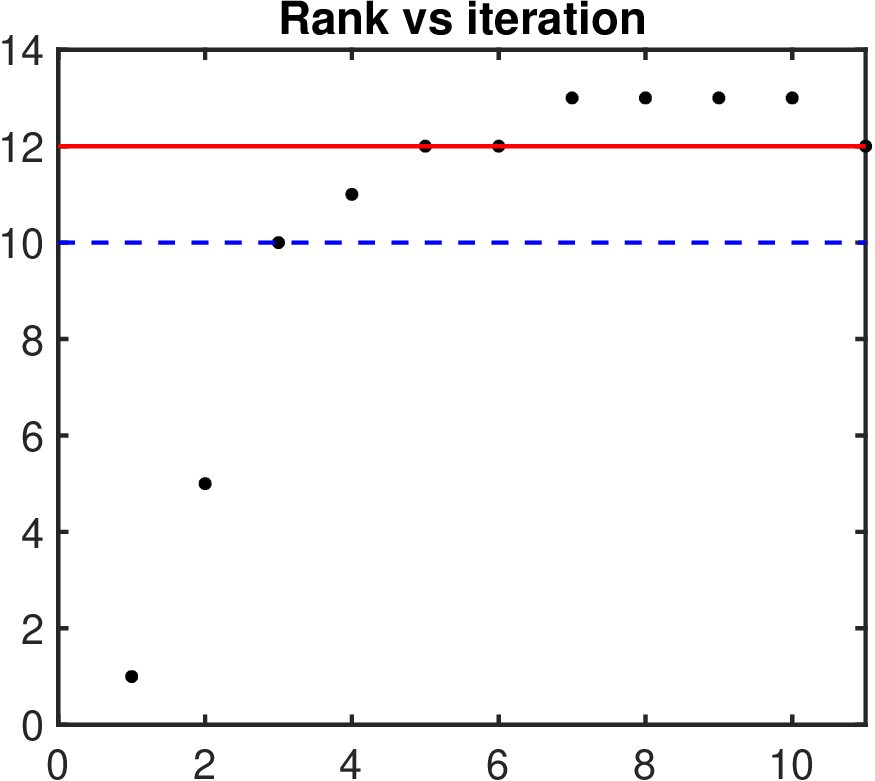}
\caption{Laplace's equation solved by ES-preconditioned lrAA. This figure displays residuals (left) and rank (right) as a function of the iteration count. The left two figures are for $n=m = 1023$ and the right two are for $n=m =4096$. Note the blue and red lines overlay each other. \label{fig:exp_sum}}
\end{center}
\end{figure}

In this example, we use much finer grids with $n=m = 1023$ and $4096$ . We use a preconditioner based on the parameters from the file \verb+Rel1_x_n10.1E10+ from the above repository. We set the rate $\alpha = 1$ in the Richardson iteration, and in this example we don't use any scheduling due to the extreme fast convergence of iteration. We terminate the iteration when the residual has been reduced by $10^{-8}$. As can be seen in Figure \ref{fig:exp_sum}, the ES-preconditioned lrAA method converges rapidly, and the intermediate rank is well-controlled.

\subsubsection{Bratu problem}
\graphicspath{{figures/bratu}}
We use lrAA to solve the non-linear Bratu problem    
\[
u_{xx} + u_{yy} + \lambda e^{u} = 0, \ \ (x,y) \in [0,1] \times [0,1],
\]
with $\lambda=1$ and homogeneous Dirichlet boundary conditions. To find the approximate solution we discretize this equation using standard second order finite difference approximations for the $x$ and $y$ derivatives. Given an approximation $X(i,j) \approx u(x_i,y_j)$ this results in a function $G_{\rm B}(i,j;X)$ describing the equation  
\begin{gather}
\begin{split}
G_{\rm B}(i,j;X) = \frac{1}{h_x^2}\left(X(i+1,j) - 2X(i,j) + X(i-1,j)  \right) \\
+ \frac{1}{h_x^2} \left(X(i,j+1) - 2X(i,j) + X(i,j-1)  \right) + \lambda e^{X(i,j)}.
\end{split}
\end{gather}
Near the boundaries some of the terms in this expression will be set to zero to account for the homogeneous Dirichlet boundary conditions. In the numerical examples below we take the mesh to be  $x_i = i h_x,$  $h_x = \frac{1}{m+1}$ and $y_j =  j h_y,$ $h_y = \frac{1}{n+1}$ with $m=n=200$, making the setup the same as  in \cite{tangAATGS}.

The fixed point function $G(i,j)$ is obtained by applying the preconditioned Richardson iteration. We have 
\[
X^{k+1}(i,j) = G(i,j;X^{k},\alpha) \equiv  X^k(i,j) + \alpha M(G_{\rm B}(i,j;X^{k})).
\] 
We test lrAA with no preconditioner and with the ES   preconditioner (corresponding to \verb+Rel1_x_n10.1E10+) described above. The lrAA parameters used are the following ${\rm TOL}=10^{-6}, \hat{m}=5, \theta=0.9, \alpha = 0.125h_x^2$ (un-preconditioned case) and   $\alpha = 0.1$ (preconditioned case, no scheduling is used). In all experiments take the initial data to be zero (represented by a rank 1 matrix).

\begin{figure}[htb]
\begin{center}
\includegraphics[width=0.3\textwidth,trim={0.0cm 0.0cm 0.0cm 0.0cm},clip]{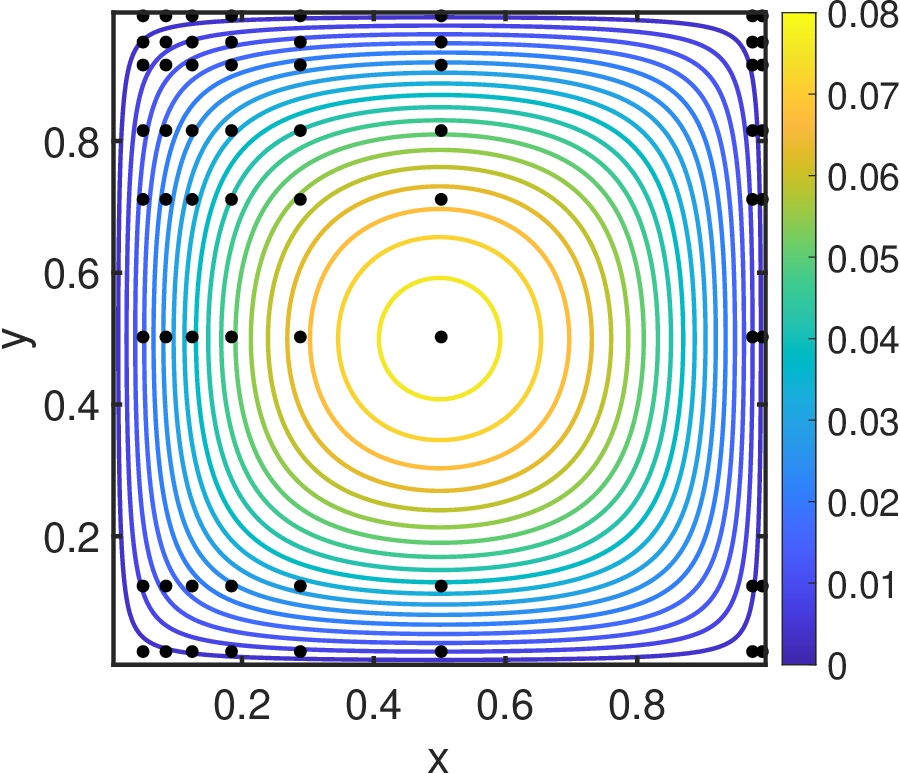} 
\includegraphics[width=0.3\textwidth,trim={0.0cm 0.0cm 0.0cm 0.0cm},clip]{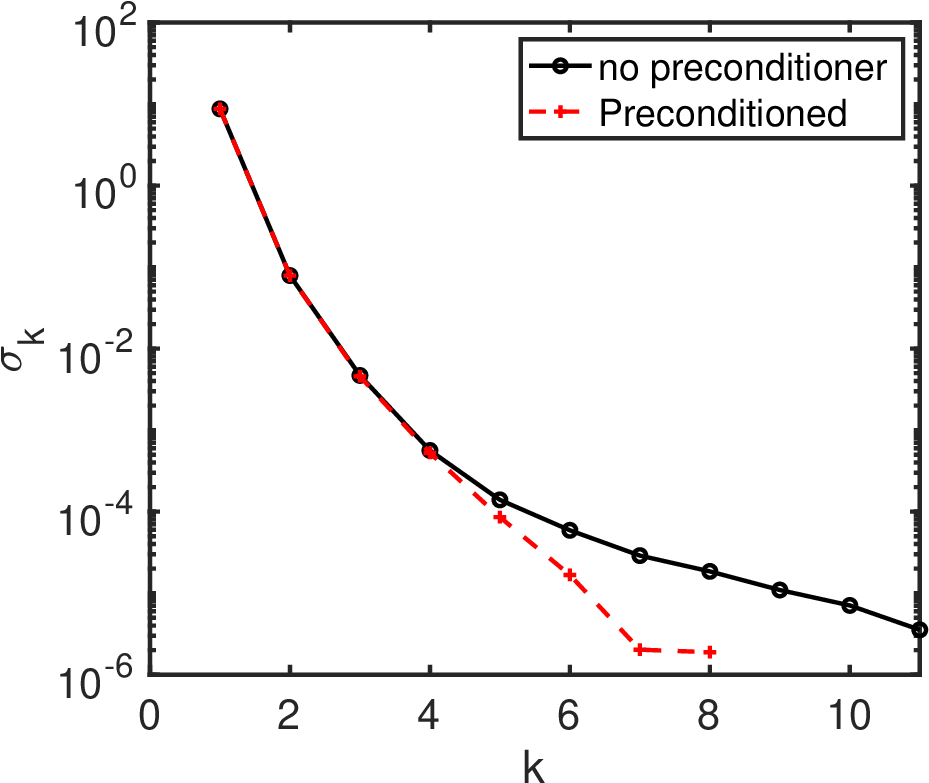}   \\ 
\includegraphics[width=0.3\textwidth,trim={0.0cm 0.0cm 0.0cm 0.0cm},clip]{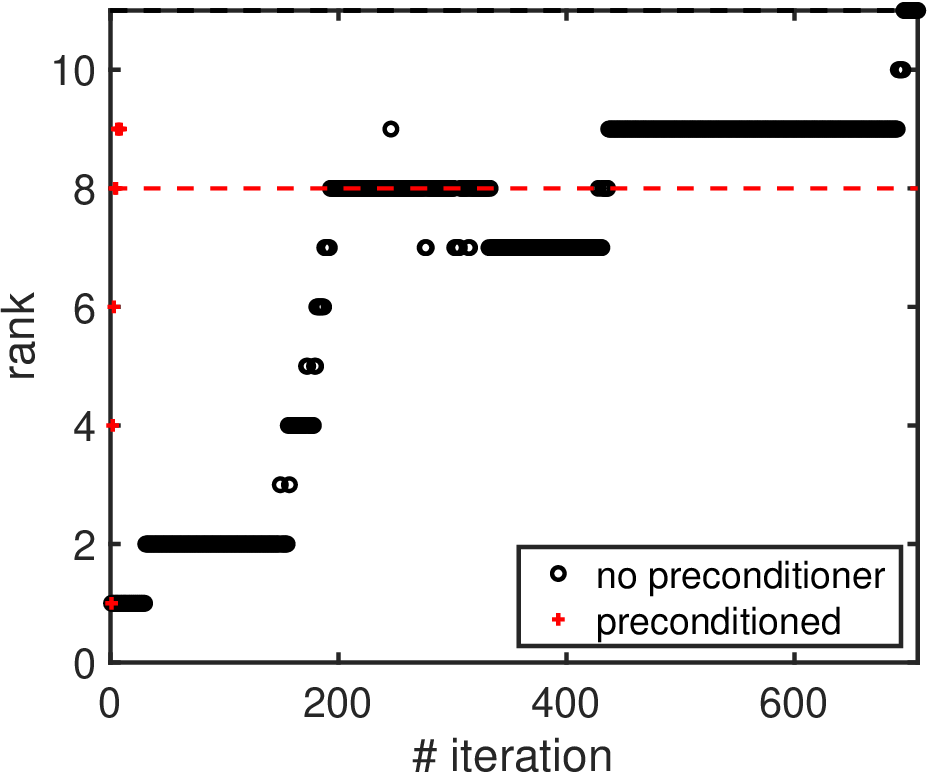} 
\includegraphics[width=0.3\textwidth,trim={0.0cm 0.0cm 0.0cm 0.0cm},clip]{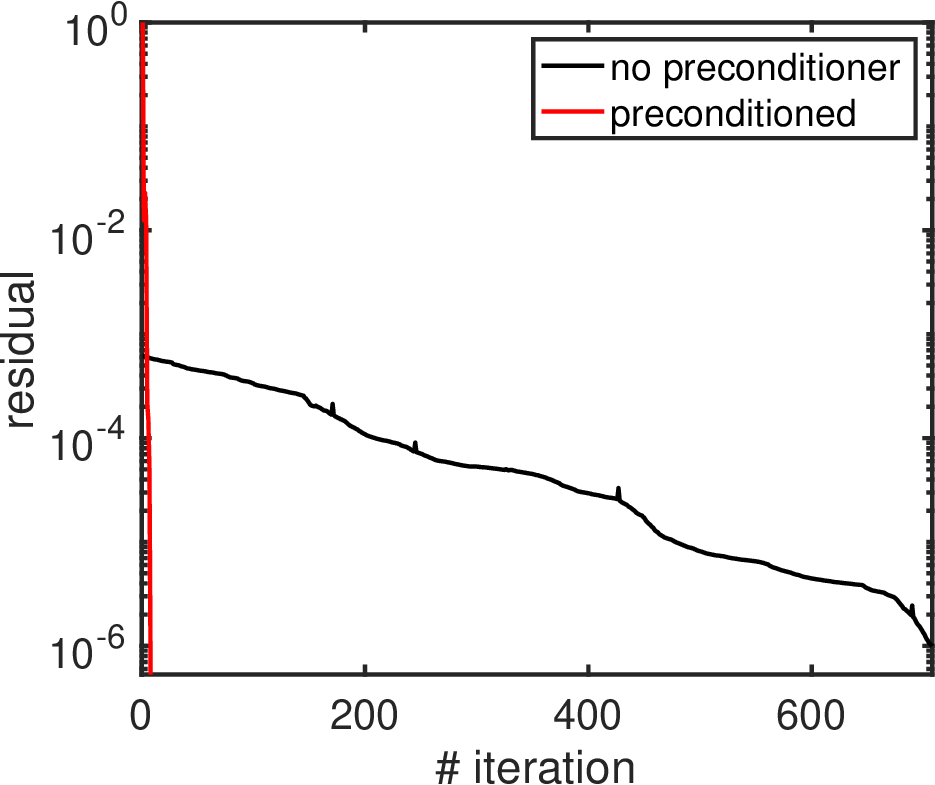}  
\caption{Bratu problem solved by lrAA. The top left figure   displays the contour levels of the converged solution along with markers at the intersection points of the final index sets $\mathcal{I}$ and $\mathcal{J}$ (for both lrAA with and without ES preconditioner). The top right figure displays the singular values from lrAA solutions with and without the ES preconditioner. The bottom left figure and right figures display  the rank evolution and the decay of the residual throughout the lrAA iterations for lrAA   with and without the ES preconditioner. 
\label{fig:bratu1}}
\end{center}
\end{figure}

The results   in Figure \ref{fig:bratu1} display the numerical solutions obtained by lrAA methods with and without the ES preconditioner. In particular, both methods obtain visually similar numerical results in terms of solution contours and column and row index section for the final iterates. 
The un-preconditioned lrAA gives a montonically increasing intermediate ranks. The ES-preconditioned lrAA converges very rapidly in 8 iterations. 

\begin{figure}[htb]
\begin{center}
\includegraphics[width=0.31\textwidth,trim={0.0cm 0.0cm 0.0cm 0.0cm},clip]{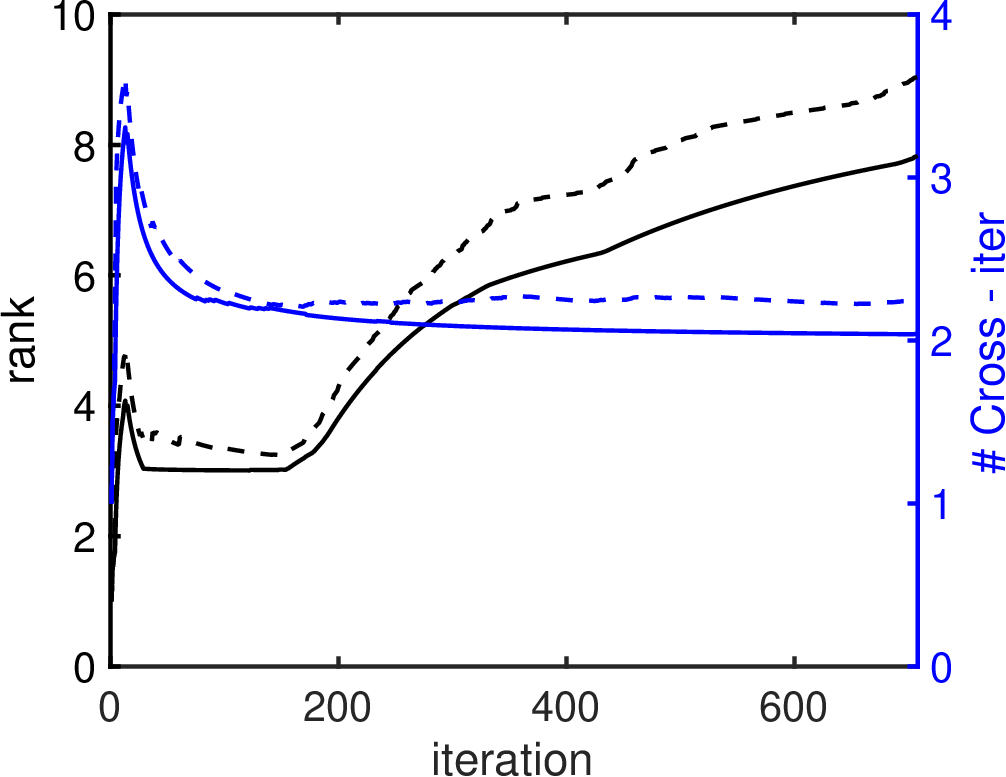}
\includegraphics[width=0.32\textwidth,trim={0.0cm 0.0cm 0.0cm 0.0cm},clip]{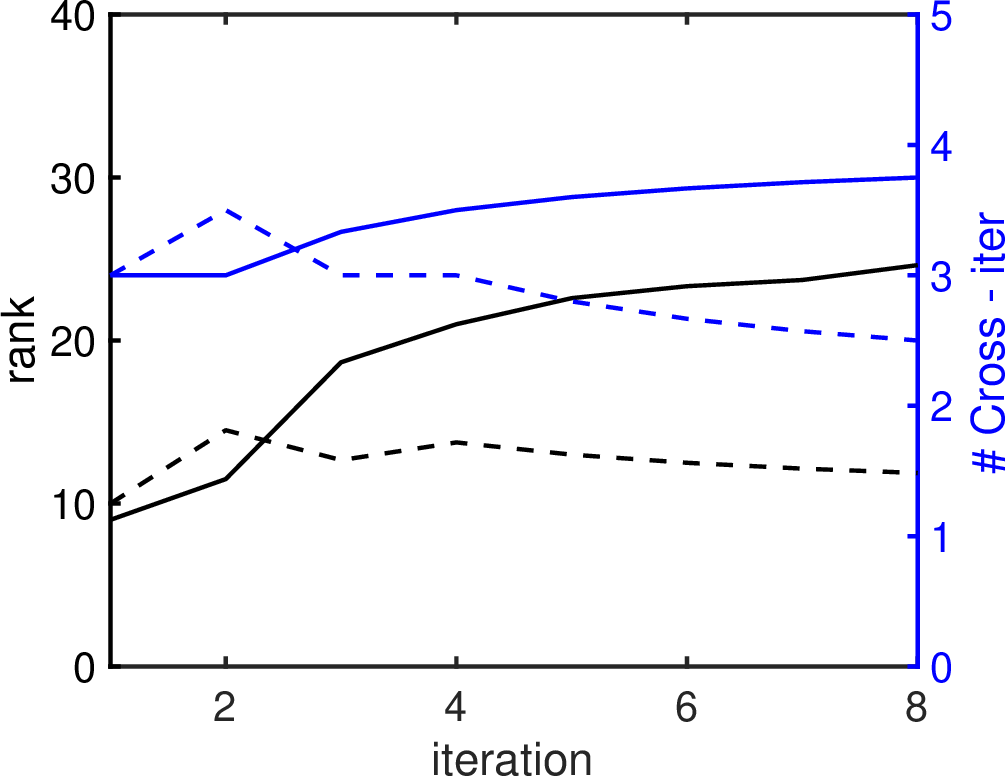}
\caption{Bratu problem solved by lrAA. Benchmarking performance of Cross-DEIM for un-preconditioned (left figure) and preconditioned (right figure) lrAA.  The figure displays the cumulative average of the max intermediate rank for Cross-DEIM (axis on the left) and the number of iterations needed for Cross-DEIM (axis on the right). The solid lines represent the results for the  fixed point function and for the dashed lines are for rounding the linear combination of previous $G$'s.  \label{fig:bratu2}}
\end{center}
\end{figure}

Next, we investigate the performance of Cross-DEIM within lrAA iterations. In particular, we plot the cumulative averages (over lrAA iterations) of the max intermediate rank and iteration number of Cross-DEIM in Figure \ref{fig:bratu2}. For both the preconditioned and un-preconditioned lrAA, we can see the iteration number of Cross-DEIM ranges from 2 to 4. The cumulative average of max intermediate rank is about 7, which is on par with the final rank (that is 10) for the un-preconditioned case. The cumulative average of max intermediate rank is about  10 for the rounding and 30 for the nonlinear evaluation, which is larger than the un-preconditioned case. However, we note that this is expected due to the larger change between the lrAA iterates due to the rapid convergence. Overall, the Cross-DEIM method works well in the lrAA methods for this example.

\subsubsection{The elliptic Monge-Amp\`{e}re equation}
\graphicspath{{figures/monge_ampere}}
As an example of the performance of lrAA for a fully nonlinear elliptic PDE, we now consider the elliptic Monge-Amp\`{e}re equation on the form (see \cite{EllipticMA})
\[
\frac{\partial^2 u}{\partial x^2}\frac{\partial^2 u}{\partial y^2} - \left(\frac{\partial^2 u}{\partial x \partial y} \right)^2  = f(x,y),
\]
on $(x,y) \in [0,1] \times [0,1]$ and with Dirichlet boundary conditions. To discretize this equation we use the scheme denoted ``Method 1'' in \cite{EllipticMA}. The scheme is defined as a fixed point iteration (here we use the traditional finite difference notation) 
\begin{equation}
u_{i,j} = h_{i,j} \equiv \frac{1}{2}(a_1+a_2) - \frac{1}{2} \sqrt{(a_1-a_2)^2 + \frac{1}{4}(a_3-a_4)^2 - h^4 f_{i,j}},
\end{equation}
where 
\[ 
2a_1 = u_{i+1,j}+u_{i-1,j}, \ \ 2a_2 = u_{i,j+1}+u_{i,j+1}, \ \ 2a_3 = u_{i+1,j+1}+u_{i-1,j-1}, \ \ 2a_4 = u_{i+1,j-1}+u_{i-1,j+1}.
\]
To iterate on this discretization we use a Richardson iteration corresponding to the fixed point iteration function
\[
G(i,j) = X(i,j) + 0.9 (H(i,j) - X(i,j)),
\]
where as before $X(i,j) = u_{i,j} \approx u(x_i,y_j)$. 

Here we use an equidistant grid with the same number of points in each direction. As in \cite{EllipticMA} we use the solution to the linear problem $u_{xx}+u_{yy} = \sqrt{2f}$ as initial data. We consider the  problem defined by the forcing 
\[
  f(x,y) = \frac{1}{\sqrt{x^2+y^2}}, 
\]
leading to the exact solution
\[
  u = \frac{2 \sqrt{2}}{3}(x^2+y^2)^{\frac{3}{4}}.
\] 
We compute the solution to this problem using lrAA on meshes with $21$, $61$, $101$ and $221$ grid points in each direction. We use a scheduling $\epsilon_{k+1} = 0.25 \rho_k$,  and take the window size $\hat{m}$ to be 5. We consider two types of tolerance for stopping the iteration, either we take the tolerance to be $10^{-10}$ or we take the tolerance to be on the order of the truncation error as measured in the Frobenius norm and set it to be a mesh-dependent constant $0.01h_x$. Here, we note that because we use a second order finite difference scheme, we expect the local truncation error to be $O(h_x^2+h_y^2),$ rescale this to the matrix Frobenious norm and with $h_x=h_y,$ it becomes $O(h_x).$ For low-rank methods for PDEs,   using a tolerance on par with the local truncation error from the underlying discretization will yield good performance in general, and avoid overly resolving the artifacts caused by the discretization schemes.

\begin{table}[htp]
\begin{center}
\begin{tabular}{|c|c|c|c|}
\hline
$n=m$ & lrAA iterations &iterations reported in\cite{EllipticMA} & final rank\\
\hline
21 &109 (7) &1083 & 13 (4) \\
61 &287 (8) &8967 & 18 (7) \\
101 &443 (13) &23849 & 20 (7)\\
221 & 675 (11) &107388 & 26 (8) \\
\hline
\end{tabular}
\caption{Iterations needed and final rank for the lrAA method for solving the elliptic Monge-Amp\`ere equation for an exact solution $u(x,y) = \frac{2 \sqrt{2}}{3}(x^2+y^2)^{\frac{3}{4}}$. The numbers corresponds to a fixed tolerance $10^{-10}$ or a local truncation error guided tolerance $0.01h_x$ (in parenthesis).\label{tab:MA}}
\end{center}
\end{table}%

\begin{figure}[htb]
\begin{center}
\includegraphics[width=0.31\textwidth,trim={0.0cm 0.0cm 0.0cm 0.0cm},clip]{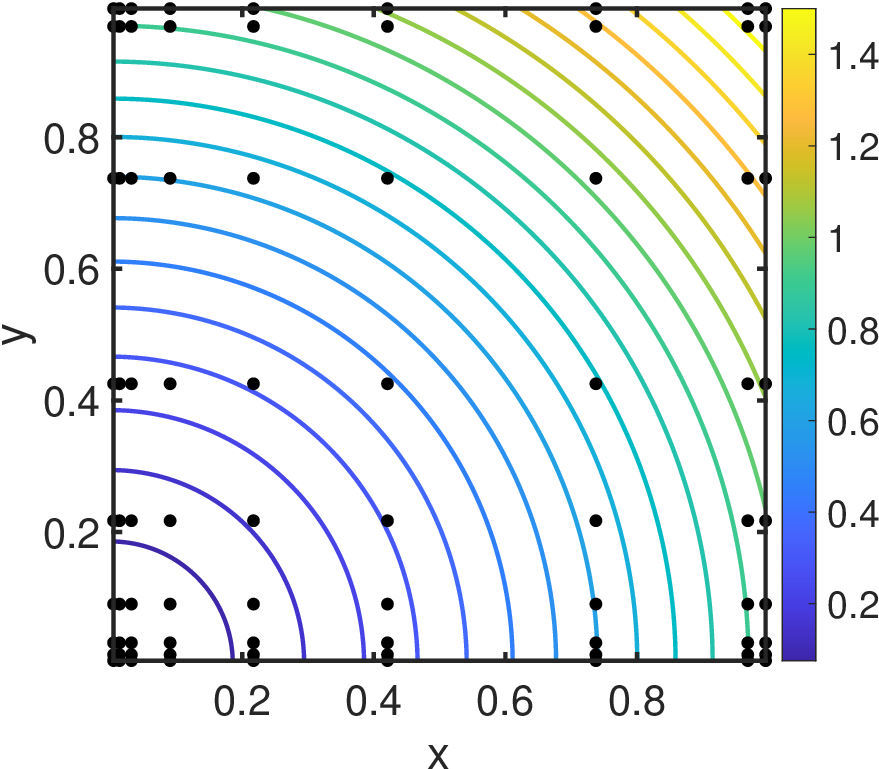}  
\includegraphics[width=0.31\textwidth,trim={0.0cm 0.0cm 0.0cm 0.0cm},clip]{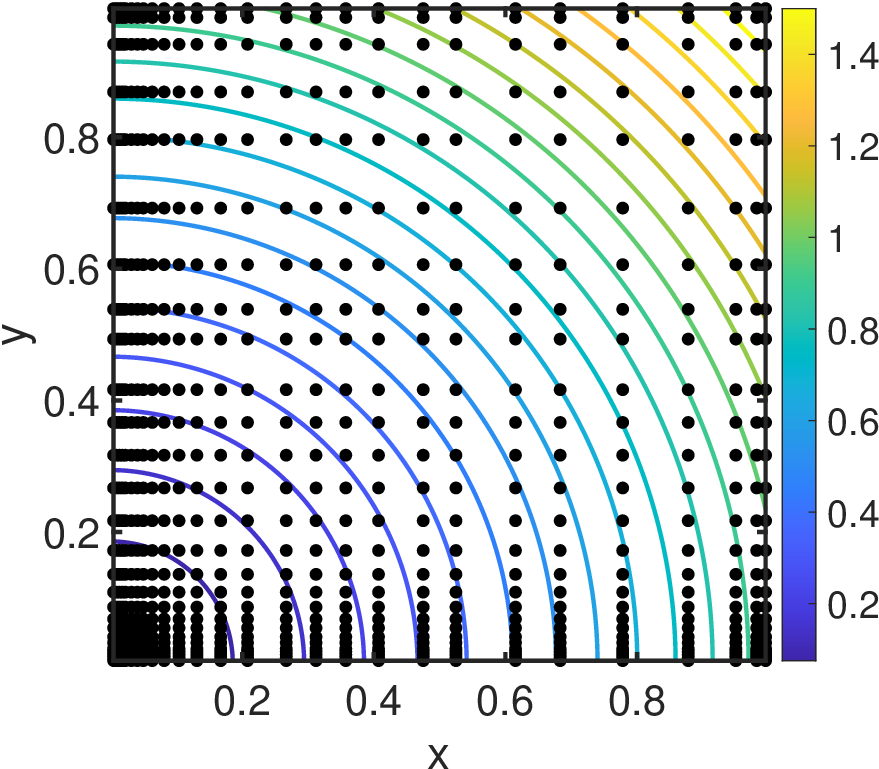} 
\caption{Elliptic Monge-Amp\`{e}re equation solved by lrAA. The results are for 221 gridpoints. Contour plot of the converged solution (for tolerance $0.01h$ (left) and $10^{-10}$ (right)) along with the index selection (black markers).   \label{fig:MA1}}
\end{center}
\end{figure}

In Table \ref{tab:MA}, we compare the number of iterations needed to reach convergence  using the Gauss-Seidel method of \cite{EllipticMA} (note that \cite{EllipticMA} uses a stricter tolerance of $10^{-14}$). The number of iterations needed for lrAA are substantially smaller than those reported in \cite{EllipticMA}.  We further note that the number of iterations and final ranks are drastically smaller when using   tolerance $0.01h_x$ that is guided by  the finite difference local truncation error.

\begin{figure}[htb]
\begin{center}
\includegraphics[width=0.31\textwidth,trim={0.0cm 0.0cm 0.0cm 0.0cm},clip]{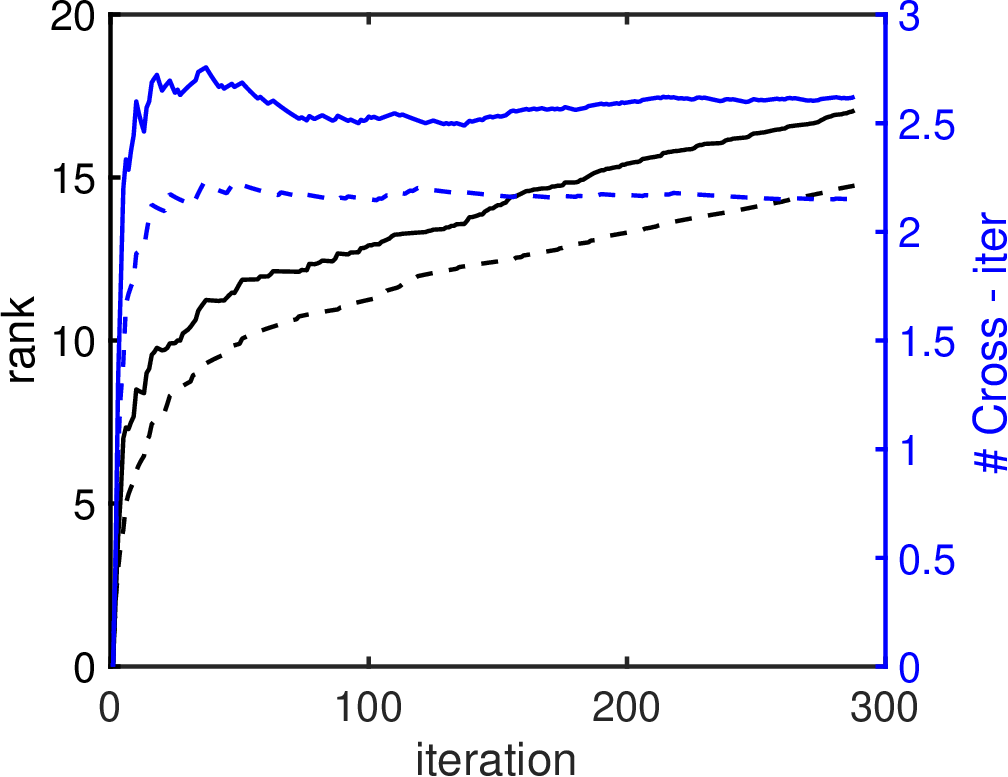} \ \
\includegraphics[width=0.31\textwidth,trim={0.0cm 0.0cm 0.0cm 0.0cm},clip]{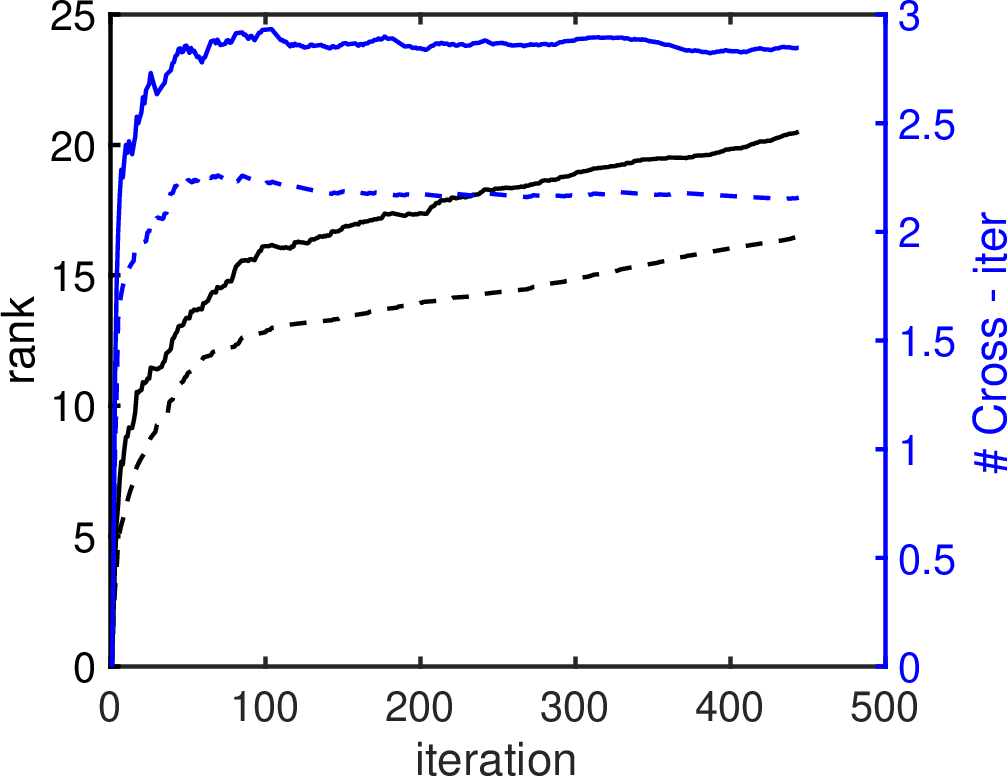} \ \ 
\includegraphics[width=0.31\textwidth,trim={0.0cm 0.0cm 0.0cm 0.0cm},clip]{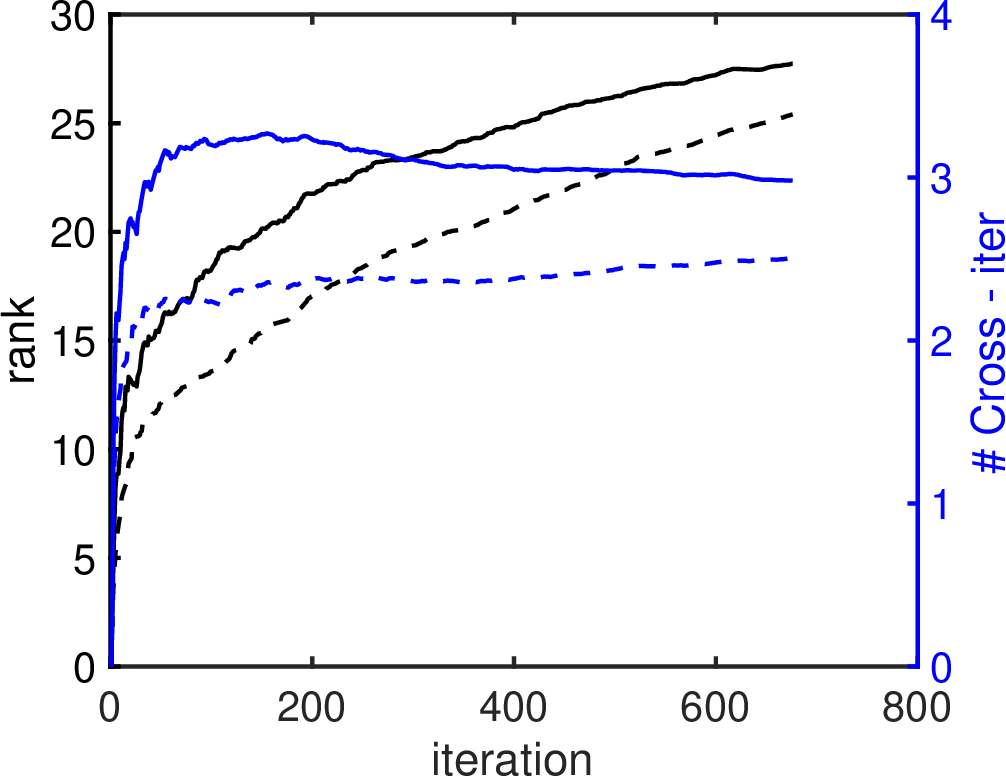}
\caption{Elliptic Monge-Amp\`{e}re equation solved by lrAA. The figure displays the cumulative average of the max intermediate rank for Cross-DEIM (axis on the left) and the number of iterations needed for Cross-DEIM (axis on the right). The solid lines represent the results for the  fixed point function and for the dashed lines are for rounding the linear combination of previous $G$'s. From left to right, the results are for 61, 101, and 221 grid points. \label{fig:MA2}}
\end{center}
\end{figure}

In Figure \ref{fig:MA1}, we display the index selection of the converged solution and the contour plots of numerical solutions for both tolerances with mesh size $n=m=221$. The results are qualitatively similar between the two tolerances.
In Figure \ref{fig:MA2}, we benchmark the performance of Cross-DEIM for this example by displaying the cumulative average of the max intermediate rank and iteration number for Cross-DEIM with various mesh size.  As can be seen, the number of iterations are very modest and the max intermediate ranks are close to the rank of the converged solution.     

\subsubsection{The Allen-Cahn equation}
\begin{figure}[htb]
\begin{center}
\includegraphics[width=0.35\textwidth,trim={0.0cm 0.0cm 0.0cm 0.0cm},clip]{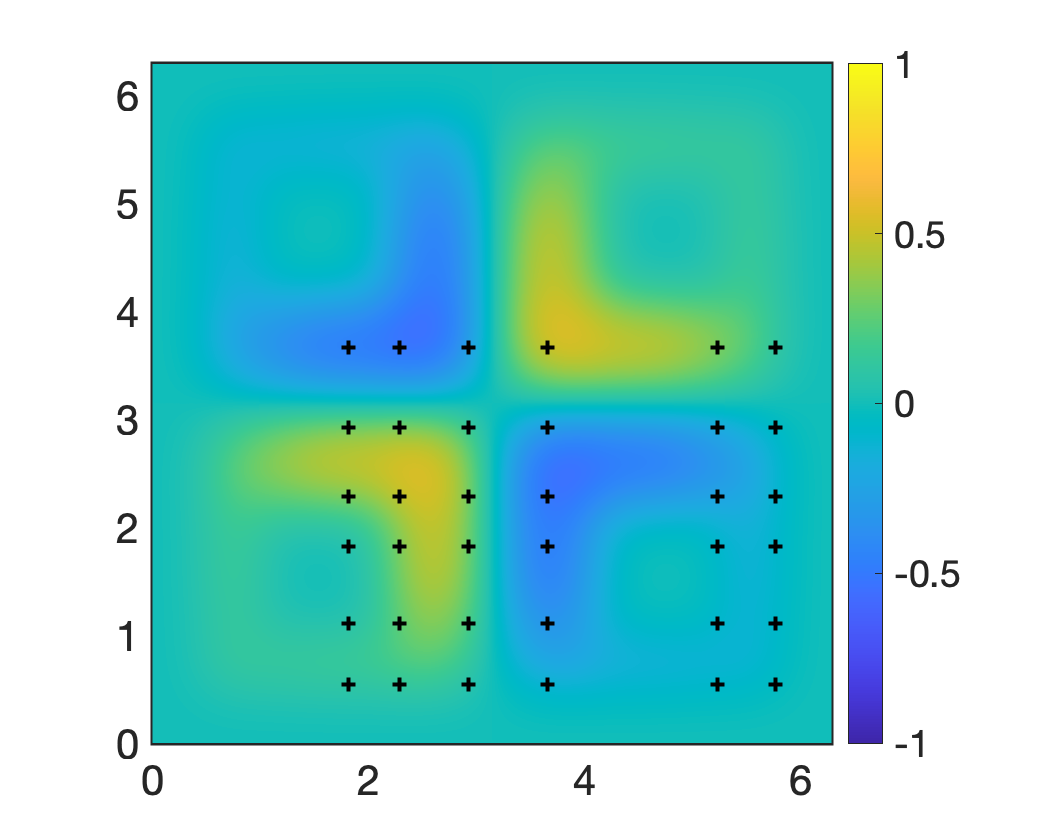} 
\includegraphics[width=0.35\textwidth,trim={0.0cm 0.0cm 0.0cm 0.0cm},clip]{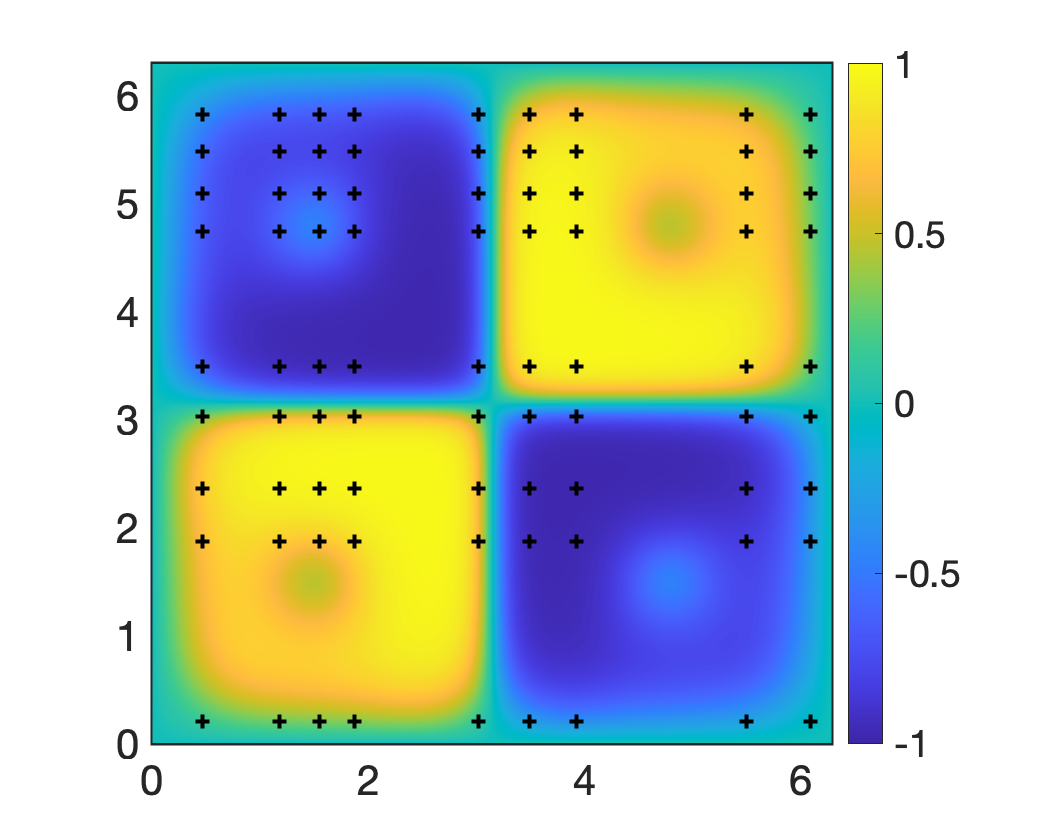} 
\includegraphics[width=0.35\textwidth,trim={0.0cm 0.0cm 0.0cm 0.0cm},clip]{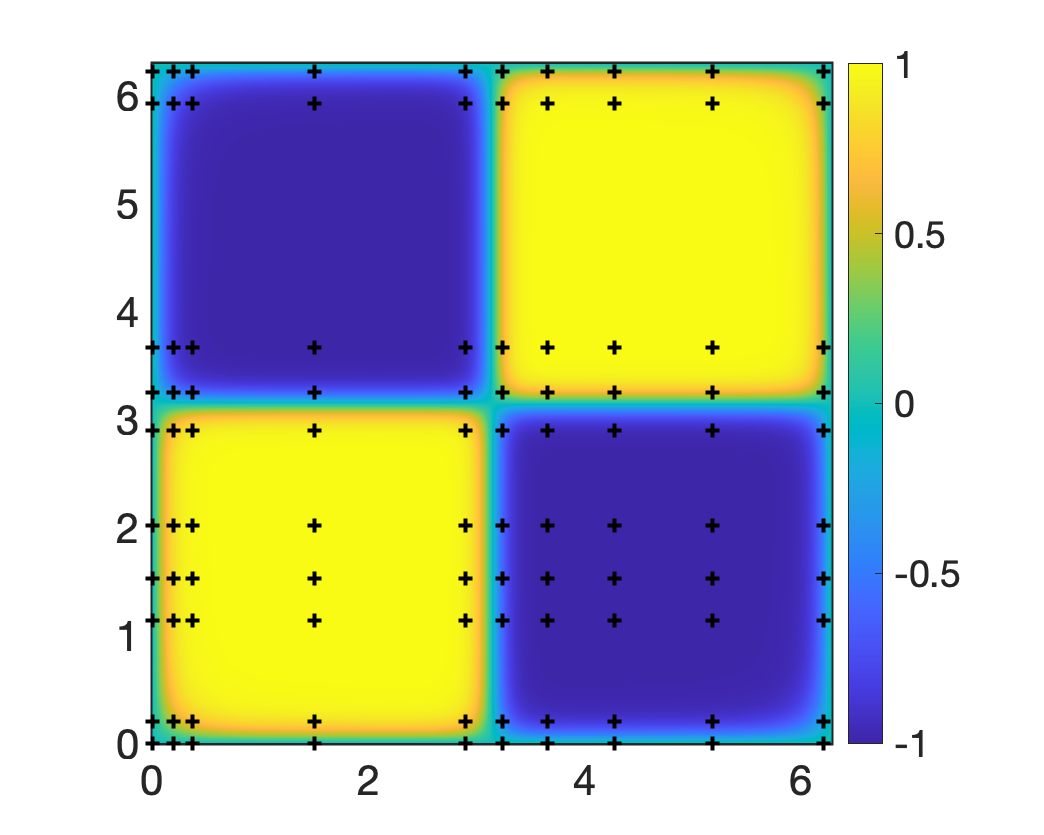} 
\includegraphics[width=0.35\textwidth,trim={0.0cm 0.0cm 0.0cm 0.0cm},clip]{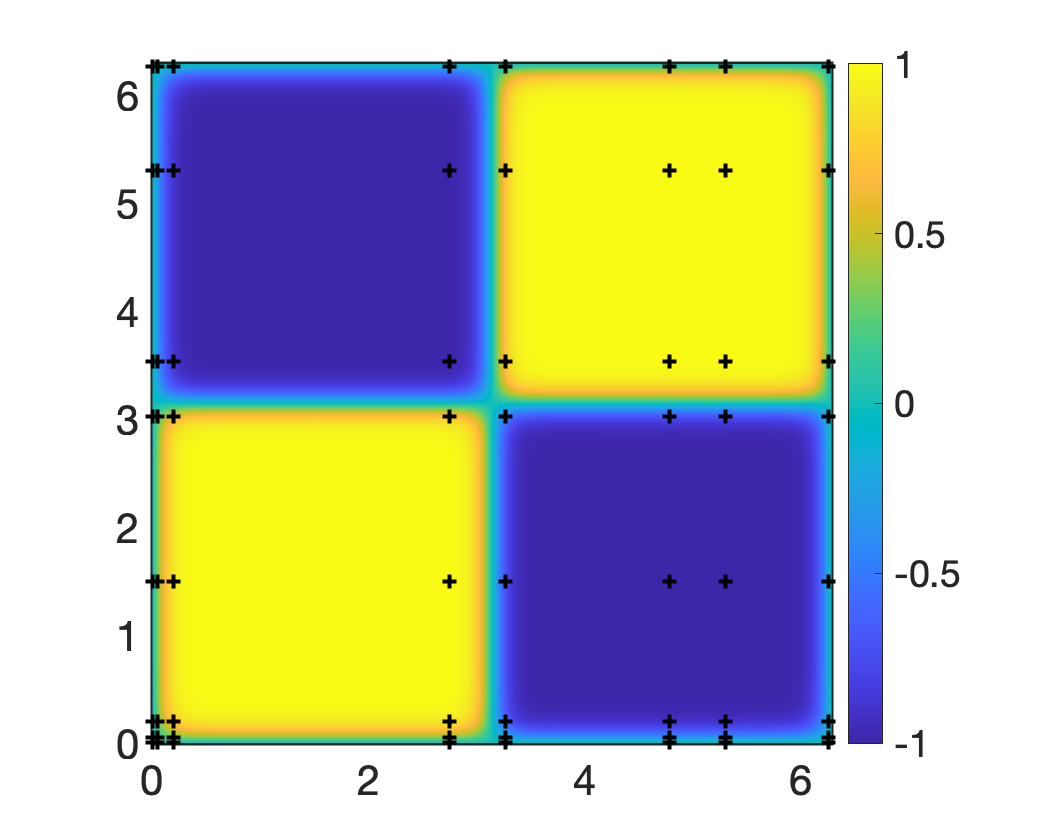} 
\caption{Numerical solution to the Allen-Cahn equation by ES-preconditioned lrAA at time 2.5, 5.0, 7.5 and 10. Markers are placed at the intersection points of the final index sets $\mathcal{I}$ and $\mathcal{J}.$ ${\rm TOL}=10^{-2}.$ \label{fig:AC1}}
\end{center}
\end{figure}

\begin{figure}[htb]
\begin{center}
\includegraphics[width=0.31\textwidth,trim={0.0cm 0.0cm 0.0cm 0.0cm},clip]{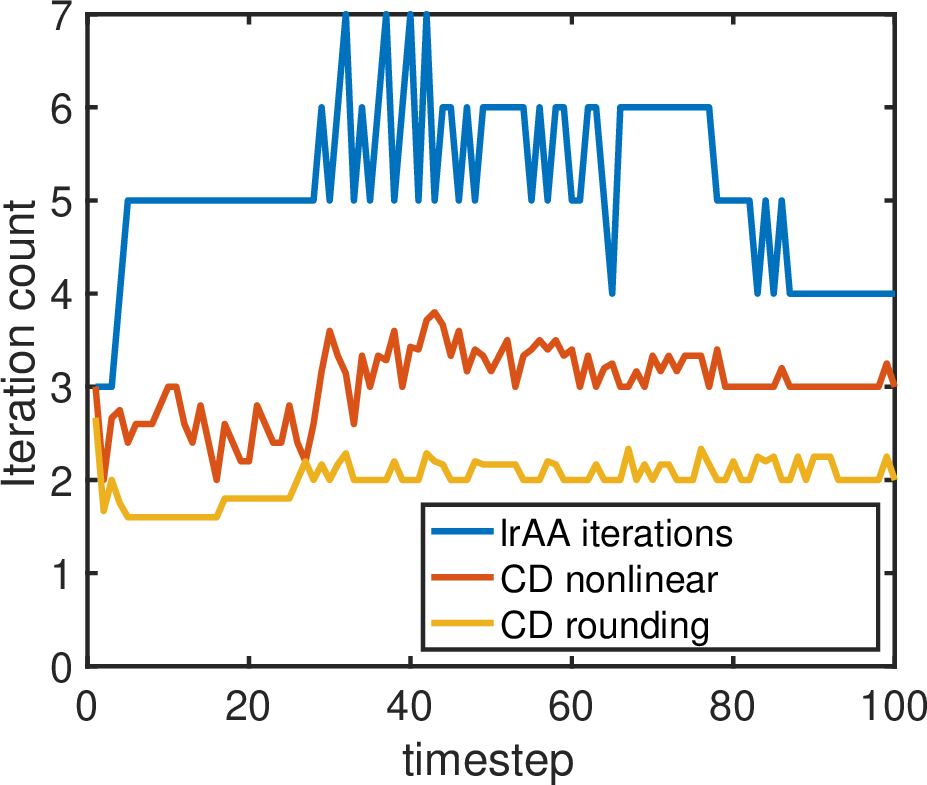} 
\includegraphics[width=0.31\textwidth,trim={0.0cm 0.0cm 0.0cm 0.0cm},clip]{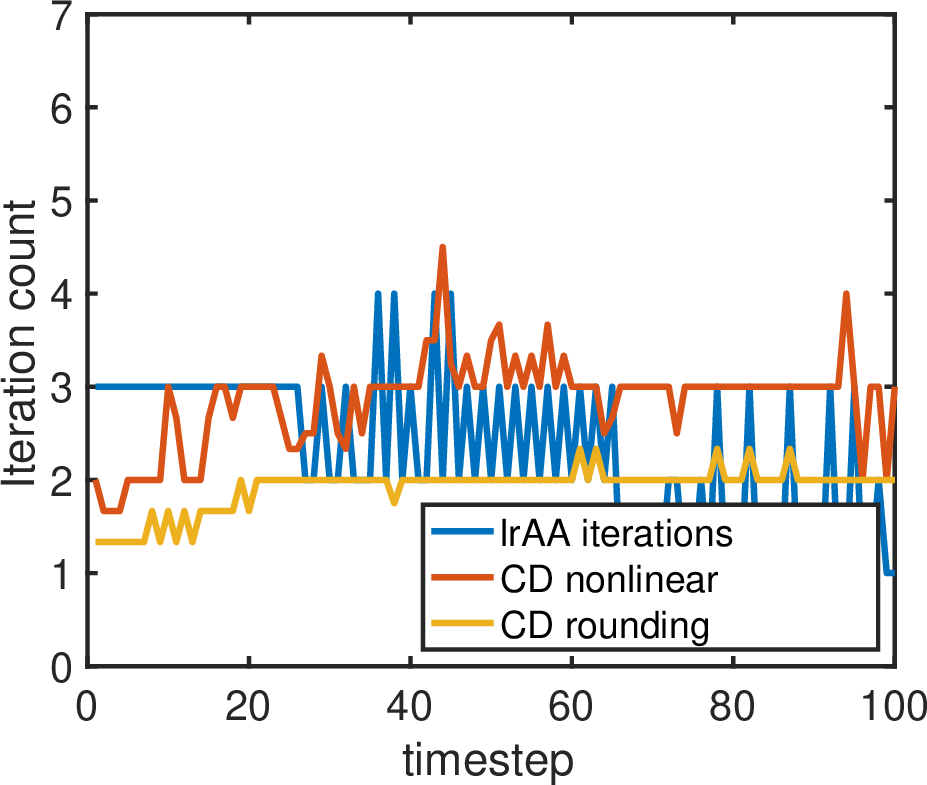} 
\caption{The Allen-Cahn equation solved by ES preconditioned-lrAA. The number of iterations used in lrAA and the two Cross-DEIM applications as a function of the time step (average values). The left is for ${\rm TOL}=10^{-4}$ and the right is for ${\rm TOL}=10^{-2}.$\label{fig:AC2}}
\end{center}
\end{figure}

Here we use the same example as in \cite{carrel2023projected,lam2024randomized} and solve the Allen-Cahn equation 
\[
u_t = \nu \Delta u +u - u^3,
\]    
with $\nu=0.01$ and periodic boundary conditions on the square domain $(x,y) \in [0,2\pi]^2$. The initial data is taken to be   
\[
    u(x,y)=\frac{[e^{-\tan^2(x)}+e^{-\tan^2(y)}]\sin(x)\sin(y)}{1+e^{|{\rm csc}(-x/2)|}+e^{|{\rm csc}(-y/2)|}}.
\]
As before $X(i,j) \approx u(x_i,y_j)$ and we again approximate the Laplacian with the standard five point stencil. We use $m = n = 256$ and solve until time 10 using 100 time steps and the backward Euler method. Here in lrAA, we take   ${\rm TOL}=10^{-2}$, in line with the local truncation error of the backward Euler method and  $\hat{m}=5$. We use the same ES preconditioner as in previous examples. In Figure \ref{fig:AC1}, we display snapshots of the solution at time 2.5, 5.0, 7.5 and 10, and the location of the Cross-DEIM points. As can be seen, they are qualitatively similar to those reported in \cite{carrel2023projected,lam2024randomized}. In Figure \ref{fig:AC2}, we display the number of iterations used in lrAA and the average number of iterations in the  two Cross-DEIM applications within lrAA for tolerances $10^{-2}$ and $10^{-4}$. As can be seen, they remain modest throughout the time stepping.

\graphicspath{{figures/Allen_Cahn}}

\section{Conclusions}
\label{sec:conclude}
In this paper, we propose lrAA, low-rank Anderson acceleration for computing low-rank solution to nonlinear matrix equations; and Cross-DEIM, an adaptive iterative cross approximation with   warm-start strategy to be used in lrAA. lrAA   only operates on the low-rank factors of the iterates, therefore   its computational cost in each iteration   scales like $O(m+n)$ instead of $O(mn).$
We propose a simple  truncation scheduling strategy that works in  controlling both the iteration number and the intermediate ranks.  The Cross-DEIM method with warm-start is demonstrated to work well  within lrAA to handle the nonlinearity.

The immediate future work is to generalize lrAA to low-rank tensor case. Another aspect is to develop more advanced scheduling strategies as those in \cite{bachmayr2017iterative} for rank truncation. We will also explore other versions of AA \cite{saad2024acceleration,tang2024anderson,lupo2024anderson} for more efficiency gains. We can also deploy Cross-DEIM for explicit low-rank schemes for time-dependent nonlinear problems.

\bibliography{lraa}
\bibliographystyle{siamplain}
\end{document}